\documentclass{article}

\usepackage[shortlabels]{enumitem}
\usepackage{amssymb}
\usepackage{amsmath}
\usepackage{amsthm}
\usepackage{amsopn}
\usepackage{graphicx}
\usepackage{bbm}
\usepackage{mathdots}
\usepackage[export]{adjustbox}
\usepackage{hyperref}
\usepackage{float}
\usepackage{nicefrac}
\usepackage{stmaryrd}
\usepackage[noadjust]{cite}

\hypersetup{
    colorlinks = true,
    citecolor=black,
    filecolor=black,
    urlcolor=black,
    linktoc=all,     
    linkcolor=blue,  
}

\graphicspath{{images/}}

\newtheorem{theorem}{Theorem}[section]
\newtheorem{definition}[theorem]{Definition}
\newtheorem{lemma}[theorem]{Lemma}
\newtheorem{remark}[theorem]{Remark}

\newtheorem{corollary}[theorem]{Corollary}
\newtheorem{proposition}[theorem]{Proposition}

\newcommand{\CC}{\mathbb C}
\newcommand{\RR}{\mathbb R}
\newcommand{\NN}{\mathbb N}
\newcommand{\ZZ}{\mathbb Z}

\newcommand{\mca}{\mathcal A}
\newcommand{\mcb}{\mathcal B}
\newcommand{\mcc}{\mathcal C}
\newcommand{\mcd}{\mathcal D}
\newcommand{\mce}{\mathcal E}

\newcommand{\mch}{\mathcal H}

\newcommand{\mct}{\mathcal T}

\newcommand{\mco}{\mathcal O}

\newcommand{\mcm}{\mathcal M}

\newcommand{\mcl}{\mathcal L}

\newcommand{\To}{\rightarrow}

\DeclareMathOperator{\re}{Re}

\DeclareMathOperator{\im}{Im}

\DeclareMathOperator{\End}{End}

\DeclareMathOperator{\Div}{div}

\DeclareMathOperator{\Id}{Id}

\DeclareMathOperator{\sub}{sub}

\DeclareMathOperator{\cl}{cl}

\DeclareMathOperator{\Vect}{Vect}
\DeclareMathOperator{\Hor}{Hor}
\DeclareMathOperator{\Ver}{Ver}
\title{On the distribution kernels of Toeplitz operators on CR manifolds} 

\author{
	Chin-Yu Hsiao\footnote{\noindent{\bf Address:} Department of Mathematics, National Taiwan University; \\ {\bf ORCID iD:} 0000-0002-1781-0013; {\bf e-mail}: chinyuhsiao@ntu.edu.tw; chinyu.hsiao@gmail.com;}
    \text{ and } Ood Shabtai\footnote{\noindent{\bf Address:} Department of Mathematics, University of Toronto; \\ {\bf ORCID iD:}  0000-0001-8970-1927; {\bf e-mail}: ood.shabtai@utoronto.ca;
}}

\begin{document}
\maketitle
\begin{abstract} We study the distribution kernel of a Toeplitz operator associated with a classical pseudodifferential operator on a compact, embeddable, strictly pseudoconvex CR manifold. The main result consists of a formula for the values at the diagonal of the second coefficient in the expansion of the symbol of the kernel. We also establish asymptotic expansions for Toeplitz operators on the positive part of a compact not necessarily strictly pseudoconvex CR orbifold under certain natural assumptions.
\end{abstract}
\tableofcontents
\section{Introduction}
Let $(X, T^{1,0}X)$ be a compact orientable strictly pseudoconvex CR manifold. The \textit{Szeg\H{o} projection} on $X$ is the orthogonal projection\footnote{Here, $L^2(X)$ is defined using a certain measure which is associated with the CR structure of $X$ (see Sect. \ref{CR_prelims}).} $$\Pi : L^2(X) \to \ker \bar \partial_b,$$
where $\bar \partial_b$ is the so-called \textit{tangential Cauchy-Riemann operator} (\ref{tcro}). Let $E$ be a classical pseudodifferential operator on $X$. Then $E$ defines a \textit{Toeplitz operator} $T_E$ by the formula \begin{equation}\label{toep_op}T_E = \Pi E \Pi.\end{equation}
The study of Toeplitz operators on CR manifolds has played a substantial role, through the seminal monograph \cite{bdmg} by L. Boutet de Monvel and V. Guillemin, in the development of the theory of quantization of K\"{a}hler manifolds (\cite{bms}). More recently (\cite{gh1, gh2, hhms}), Toeplitz operators on CR manifolds have been examined with the aim of formulating an analogous theory of quantization for CR manifolds. This aim also largely motivates the present work (see Sect. \ref{motiv_sect}), which is devoted to the study of the distribution kernels of Toeplitz operators such as $T_E$ (\ref{toep_op}). 

The distribution kernel of $\Pi$, termed the \textit{Szeg\H{o} kernel}, is the fundamental example of a distribution kernel of a Toeplitz operator. Accordingly, it will be instructive to recall some of its properties (see Sect. \ref{CR_prelims} for further details). Let $\omega_0 \in \Omega^1(X)$ be a contact form corresponding\footnote{That is, $\omega_0$ satisfies $(\ker \omega_0) \otimes \CC = T^{1,0} X \oplus \overline{T^{1,0}X}$, and $d\omega_0 \big{|}_{\ker \omega_0}$ is non-degenerate.} to the CR structure of $X$. Let $\mct$ be the Reeb vector field determined by the equations $\iota_{\mct} d\omega_0 = 0$ and $\omega_0(\mct) = -1$. 
Recall that a compact, strongly pseudoconvex CR manifold $X$ is embeddable\footnote{i.e., there exists a CR embedding of $X$ into $\CC^{N}$ for some $N$.} if and only if the range of the \textit{Kohn Laplacian} $\square_b = \bar \partial_b^* \bar \partial_b$
is closed in $L^2(X)$ (\cite{bdm, kohn1, kohn2}). For $D \subset X$ open and $f\in C^\infty (D\times D)$, we denote the differentials of $f$ with respect to its first and second arguments by $d_x f$ and $d_y f$, respectively, and similarly for other differential operators on $D$ (for example, $\mct_xf$ and $\mct_yf$). The following description of the Szeg\H{o} kernel $\Pi(x,y)$ was established in the classic work \cite{bdms} of J. Sj\"{o}strand and L. Boutet de Monvel (see also \cite{hsiao1, hsiaomari2}).
\begin{theorem}[\cite{bdms}]\label{bdms_thm} Let $X$ be an orientable compact strongly pseudoconvex embeddable CR manifold of dimension $2n+1$, $n \ge 1$. Assume that $\CC TX$ is equipped with the \textit{Levi metric} (\ref{levi_metric}). Let $(D,x)$ be a coordinate patch on $X$. Then
\begin{equation}\label{szego_ker} \Pi(x,y) \equiv \int_0^\infty e^{i t \phi(x,y)} a(x,y,t)dt \mod C^\infty(D \times D),\end{equation}
where $\phi \in C^\infty(D \times D)$ satisfies
\begin{equation*}\label{CR_flat_diag}\mbox{$\overline{\partial}_{b,x}\phi$ and 
$\partial_{b,y}\phi$ vanish to infinite order at $x=y$},\end{equation*}
and
\begin{equation}\label{phi_props} \begin{aligned} 
& \im(\phi) \ge 0,\\ & \phi(x,y) = 0 \ \text{if and only if }x=y,\\ & d_x\phi(x,x) = - d_y \phi(x,x) = -\omega_0(x),\end{aligned} \end{equation}
and $a\in S^n_{\cl}(D \times D \times \RR_+)$ is a \textit{classical H\"{o}rmander symbol} (see Definition \ref{clhs}),
\begin{equation}\label{a_props}
\begin{aligned} &a(x,y,t) \sim \sum_{j=0}^\infty a_j(x,y)t^{n-j}\ \text{in } S^n_{1,0}(D \times D \times \RR_+),\\ & a_0(x,x) = \frac 1 {2\pi^{n+1}}  \ \text{for every }x \in D.\end{aligned}\end{equation}\end{theorem}
The phase function $\phi$ and the symbol $a$ in the description of $\Pi(x,y)$ as an oscillatory integral are generally not determined uniquely. At the same time, various requirements may be imposed on $\phi$ or on $a$ in order to achieve respective uniqueness properties (cf. Lemma \ref{pointwise_uniqueness}). Similarly, by imposing further assumptions on $\phi$ it is possible to ensure that particular choices of $a$ can be made. This type of reasoning was applied in \cite{hsiao_shen} to obtain a more refined description of $\Pi(x,y)$, as follows.
\begin{theorem}[\cite{hsiao_shen}]\label{hsiao_shen_thm_intro} Let $X$ be a compact orientable strongly pseudoconvex embeddable CR manifold of dimension $2n+1$, $n \ge 1$. Let $(D,x)$ be a coordinate patch on $X$. Assume that $\hat \phi \in C^\infty(D \times D)$ satisfies (\ref{phi_props}), and $\mct_y^2 \hat \phi(x,x) = 0$ for all $x \in D$, and that $\hat \phi$ is equivalent, in the sense of \cite[Sect. 4]{ms}, to the phase function $\phi$ of Theorem \ref{bdms_thm}. Then there exists $\hat a \in S^n_{\cl}(D \times D \times \RR_+)$ such that $$\Pi(x,y) \equiv \int_0^\infty e^{it \hat \phi(x,y)} \hat a(x,y,t)dt \mod C^\infty(D \times D),$$
and $$\begin{aligned} &\hat a \sim \sum_{j \ge 0}^\infty \hat a_j t^{n-j}\ \text{in } S^n_{1,0}(D \times D \times \RR_+),\\ &\hat a_0(x,x) = \frac 1 {2\pi^{n+1}},\ \mct_y \hat a_0(x,x) = 0\text{ on }D,\\ &\hat a_1(x,x) = \frac{R(x)}{4\pi^{n+1}}. \end{aligned}$$
Here, $R$ denotes the \textit{Tanaka-Webster scalar curvature} (\ref{tw_curv}) associated with $\omega_0$.\end{theorem}
A phase function satisfying the assumptions of Theorem \ref{hsiao_shen_thm_intro} may be obtained as follows. Let $(D,x)$ be a coordinate patch such that $\mct = -\partial_{x_{2n+1}}$, and assume that $\phi \in C^\infty(D \times D)$ is the phase function of Theorem \ref{bdms_thm}. By the Malgrange preparation theorem (\cite[Theorem 7.5.5]{hormander1}), there exist $f, g \in C^\infty(D \times D)$ such that $f(x,x) = 1$ for all $x \in D$ and $$\phi(x,y) = f(x,y)(y_{2n+1} + g(x,y')),\ y' = (y_1, ..., y_{2n}).$$
Then $\hat \phi(x,y) = y_{2n+1} + g(x,y')$ satisfies (\ref{phi_props}), $\mct_y^2 \hat \phi = 0$, and $\hat \phi$, $\phi$ are equivalent in the sense of \cite{ms}. Consequently, Theorem \ref{hsiao_shen_thm_intro} applies to the phase function $\hat \phi$.
\subsection{Main results}
The main result of the current paper (Theorem \ref{main_thm}) is an extension of the description of the Szeg\H{o} kernel provided in Theorem \ref{hsiao_shen_thm_intro} to kernels of Toeplitz operators $T_E$, where $E$ is a classical pseudodifferential operator. We note that our work substantially relies on, and may be considered as a continuation of, the work carried out in \cite{hsiao_shen}.

Let us introduce several notions that will be used in the formulation of Theorem \ref{main_thm}. We refer the reader to Sect. \ref{prelim_sect} for a more detailed exposition. Let $E \in L^m_{\cl}(X)$ be a classical pseudodifferential operator of order $m \in \RR$ with principal symbol $e_0 \in C^\infty(T^*X \setminus 0)$ (here $0 \subset T^*X$ denotes the zero section). In addition to $e_0$, it is possible to associate to $E$ a type of subprincipal symbol $e_{\sub} \in C^\infty(T^* X \setminus 0)$, as follows\footnote{See Lemma \ref{density_subprinc_symb} for further details.}. Let $(D,x)$ be a coordinate patch on $X$ and let $(x, \xi)$ be the associated coordinates on $T^*D$. Assume that the volume form $dv_X$ induced from the \textit{Levi metric} (\ref{levi_metric}) is specified on $D$ by $$dv_X = \lambda(x) dx_1 \wedge ... \wedge dx_{2n+1},\ \lambda \in C^\infty(D).$$ Assume that the total symbol $e$ of $E$ over $D$ has the expansion $$e(x,\xi) \sim \sum_{j \ge 0} e_j(x,\xi),$$
where $e_j$ is homogeneous of order $m-j$ in the $\xi$ variables, $|\xi| \ge 1$. Then $e_{\sub}$ is specified on $T^*D\setminus 0$ by $$e_{\sub}(x,\xi) = e_1(x,\xi) + \frac i 2 \sum_{j=1}^{2n+1} \partial_{x_j, \xi_j} e_0(x,\xi) +\frac i {2}\sum_{j=1}^{2n+1} \partial_{\xi_j} e_0(x,\xi) \partial_{x_j}(\log |\lambda|)(x).$$
Notably, the definition of $e_{\sub}$ depends on the choice of $\omega_0 \in \Omega^1(X)$.

Next, define a differential operator $P : C^\infty(T^*X) \to C^\infty(T^*X)$, as follows\footnote{See Sect. \ref{diffr_op_def} and Corollary \ref{div_jhamvf} for further details.}. The contact form $\omega_0$ gives rise to the \textit{Tanaka-Webster connection} (see Definition~\ref{tan_web_con}). As with any affine connection, we obtain a splitting of the tangent bundle of the cotangent
bundle into horizontal and vertical subbundles, $$T T^*X = \Hor(T^*X) \oplus \Ver(T^*X) \simeq \pi^*_{T^*X}TX \oplus \pi^*_{T^*X}T^*X.$$ Here, $\pi_{T^*X} : T^*X \to X$ is the standard projection.
Now, the \textit{Levi distribution} $HX \subset TX$ is the unique subbundle of $TX$ such that \begin{equation}\label{levi_distribution}HX \otimes \CC = T^{1,0}X \oplus T^{0,1}X,\ T^{0,1}X = \overline{T^{1,0}X}.\end{equation} The Levi distribution is equipped with a complex structure (\ref{cplx_str}) $J : HX \to HX$. Let $H^*X \subset T^*X$ be the dual of $HX$, equipped with the dual complex structure $J : H^*X \to H^*X$. Let $H X_{\Hor} \subset \Hor(T^*X)$ and $J_{\Hor} : H X_{\Hor} \to H X_{\Hor}$ denote the horizontal lifts of $HX$ and $J:HX \to HX$, respectively. Similarly, let $H X_{\Ver} \subset \Ver(T^*X)$ and $J_{\Ver} : HX_{\Ver} \to HX_{\Ver}$ denote the vertical lifts of $H^*X$ and $J : H^*X \to H^*X$. Let $\mct^{\Hor}$ and $\omega_0^{\Ver}$ denote the horizontal lift of $\mct$ and the vertical lift of $\omega_0$, respectively. Define a morphism $J_{T^*X} : T T^*X \to T T^*X$ by \begin{equation*}\begin{aligned} &J_{T^*X} \big{|}_{HX_{\Hor}} = J_{\Hor},\\ &J_{T^*X} \big{|}_{HX_{\Ver}} = - J_{\Ver},\\ &J_{T^*X} \mct^{\Hor} = J_{T^*X} \omega_0^{\Ver} = 0. \end{aligned}\end{equation*}
Next, recall that $T^*X$ is equipped with a natural symplectic form.
Let $X_F$ denote the Hamiltonian vector field of a function $F \in C^\infty(T^*X)$. Finally, set $$P(F) = -\frac 1 2\Div(J_{T^*X} X_F),$$
where $\Div : \Vect(T^*X) \to C^\infty(T^*X)$ denotes the divergence with respect to the Liouville volume form.

The main result of the present article is as follows.
\begin{theorem}\label{main_thm} Let $(X, T^{1,0}X)$ be an orientable compact embeddable strongly pseudoconvex CR manifold of dimension $2n+1$, $n\ge 1$. Let $\omega_0 \in \Omega^1(X)$ be a contact form, as specified in (\ref{contact_form}), and let $dv_X$ be the associated volume form\footnote{See Sect. \ref{CR_prelims}}. Let $(D,x)$ be a coordinate patch on $X$, and let $(x,\xi)$ denote the induced coordinates on $T^*D$. Assume that $\hat \phi \in C^\infty(D \times D)$ satisfies (\ref{phi_props}), $\mct_y^2 \hat \phi(x,x) = 0$ for all $x \in D$, and $\hat \phi$ is equivalent, in the sense of \cite{ms}, to the phase function $\phi$ of Theorem \ref{bdms_thm}. Let $m \in \RR$ and $E \in L^m_{\cl}(X)$. Define $\mce_0(x) = e_0(x, -\omega_0(x))$, where $e_0 \in C^\infty(T^*X \setminus 0)$ is the principal symbol of $E$.  Let $T_E(x,y)$ denote the distribution kernel of the Toeplitz operator $T_E = \Pi E \Pi$. Then
\begin{equation*} T_E(x,y) \equiv \int_0^\infty e^{it  \hat \phi(x,y)} \hat b(x,y,t) dt,\end{equation*}
where $\hat b(x,y,t) \in S^{n+m}_{\cl}(D \times D \times \RR_+)$ satisfies
\begin{equation*} \hat b(x,y,t) \sim \sum_{j=0}^\infty \hat b_{j}(x,y) t^{n+m-j}\text{ in } S^{n+m}_{1,0}(D \times D \times \RR_+)\end{equation*}
with \begin{equation}\label{main_formula} \begin{aligned} \hat b_{0}(x,x) &= \frac{e_0(x,-\omega_0(x))}{2\pi^{n+1}},\ \mct_y \hat b_0(x,x) = 0\text{ on }D,\\ \hat b_{1}(x,x) &= \frac 1 {4\pi^{n+1}}\left[R(x) \mathcal E_0(x) -\square_b \mathcal E_0(x)  + P(e_0)(x, -\omega_0(x))\right]\\&+\frac 1 {4\pi^{n+1}} \left[2 e_{\sub}(x, -\omega_0(x))-i m \mct \mce_0(x) \right].\end{aligned} \end{equation} Here, $P :C^\infty(T^*X) \to C^\infty(T^*X)$ is the differential operator specified in Corollary \ref{div_jhamvf} and $e_{\sub} \in C^\infty(T^*X\setminus 0)$ is the subprincipal symbol of $E$ with respect to the $1$-density $|dv_X|$, as specified in Lemma \ref{density_subprinc_symb}.\end{theorem}
\begin{remark} A key feature of Theorems \ref{hsiao_shen_thm_intro}, \ref{main_thm} is that regardless of the use of local coordinates, the coefficients in the expansions of the symbols of $\Pi$, $T_E$ are expressed in terms of globally defined objects. In fact, it is possible to formulate Theorem \ref{main_thm} using a globally defined phase function and a globally defined symbol.\end{remark}
We note that the representation of $T_E(x,y)$ formulated in Theorem \ref{main_thm} satisfies certain limited uniqueness properties, as specified in Corollary \ref{uniqueness_corollary} (cf. \cite[Lemma 1.1]{hsiao_shen}).
Notably, if $E = \mcm_f$ is the operator of multiplication by a smooth function $f \in C^\infty(X)$, then formula (\ref{main_formula}) simplifies considerably.
\begin{corollary} Let $E = \mcm_f$, where $f \in C^\infty(X)$. Then using the notations of Theorem \ref{main_thm}, we have that $m = 0$ and $e_0(x,\xi) = f(x) = \mce_0(x)$, while\footnote{Note that $X_f \in \Ver(T^*X)$, hence $P(f) = 0$.} $$e_{\sub} = 0 = P(f).$$ Hence, $$\begin{aligned} &\hat b_0(x,x) = \frac{f(x)}{2\pi^{n+1}},\\ & \hat b_1(x,x) = \frac 1 {4\pi^{n+1}} \left[R(x) f(x) - \square_b f(x) \right]. \end{aligned}$$\end{corollary}

In Section~\ref{s-gue251129yyd}, we establish asymptotic expansions for Toeplitz operators on the positive part of a compact not necessarily strictly pseudoconvex CR orbifold under certain natural assumptions (see Theorem~\ref{t-gue251209yydf}).
\subsection{Deformation quantization of CR manifolds}\label{motiv_sect}
The work presented in this paper is partly motivated by questions arising in the context of quantization of CR manifolds (\cite{fitzpatrick, hmm, gh1, gh2, hhms}), a topic which is closely related to that of quantization of symplectic manifolds, and K\"{a}hler manifolds in particular (\cite{bms, charles, charles1, mamarinescu, lefloch, schlichenmaier}).
Consider a closed connected K\"{a}hler manifold $(M, \omega)$ of complex dimension $n$. Assume that $L \to M$ is a holomorphic Hermitian line bundle such that the curvature of its Chern connection equals $-i\omega$. Let $\mch_k$ denote the space of holomorphic sections of $L^{\otimes k}$. Then $\mch_k$ is equipped with a natural inner product, obtained by integrating the Hermitian product of $L$ with respect to the Liouville measure\footnote{That is, the measure $\frac{|\omega^{\wedge n}|}{n!}$}. In geometric quantization, $\mch_k$ is viewed as the space of states of a quantum system which ``corresponds'' to $M$, and $k \to \infty$ is the semiclassical limit. The quantum counterpart of a classical observable $f \in C^\infty(M)$ is a type of Toeplitz operator $T_k(f) \in \End(\mch_k)$. The quantization maps $$T_k : C^\infty(M) \to \End(\mch_k)$$ satisfy several interesting properties. In particular, for any $f, g \in C^\infty(M)$, $$T_k(f) T_k(g) \sim \sum_{j \ge 0} T_k(C_j(f,g)) \frac 1 {k^j}$$
as $k \to \infty$, where $C_j$, $j \ge 0$, is a bidifferential operator, \begin{equation}\label{star_prod_prop}\begin{aligned} &C_0(f,g) = fg\\ &C_1(f,g) - C_1(g,f) = i \{f,g\}.\end{aligned}\end{equation}Here, $\{f,g\}$ denotes the Poisson bracket of $f, g \in C^\infty(M)$.

The quantization of $(M,\omega)$ can also be described through the language of CR geometry. The unit circle bundle $X$ in the dual bundle of $L$ is a strictly pseudoconvex CR manifold. The Hilbert spaces $\mch_k$ may be identified with spaces of ($k$-th degree) equivariant CR functions on $X$. The properties (\ref{star_prod_prop}) were established in \cite{bms} through the study of Toeplitz operators on $X$, using the theory developed in \cite{bdmg}. The properties (\ref{star_prod_prop}) express the fact that the product given by \begin{equation}\label{star_prod}f*g = \sum_{j \ge 0} C_j(f,g)\hbar^j\end{equation}
is a star product on $C^\infty(M)[[\hbar]]$ (the algebra of formal power series in $\hbar$ with coefficients in $C^\infty(M)$), so $\left(C^\infty(M)[[\hbar]], *\right)$ is a deformation quantization. Recall that by definition, a star product on $C^\infty(M)[[\hbar]]$ is an associative, unital product of the form (\ref{star_prod}), such that $C_0$, $C_1$ satisfy (\ref{star_prod_prop}).

The examples presented by geometric and deformation quantization of K\"{a}hler manifolds motivate the study of similar topics in the setting of CR manifolds. Let $(X, T^{1,0}X)$ be a compact orientable strictly pseudoconvex CR manifold. Then for any\footnote{Instead of $f$, $g$, we could also take pseudodifferential operators $F \in L^{m_1}_{\cl}(X)$, $G \in L^{m_2}_{\cl}(X)$.} $f, g \in C^\infty(X)$, the Toeplitz operators\footnote{For $f \in C^\infty(X)$ we abbreviate $T_f := T_{\mcm_f}$.} $T_f$, $T_g$ may be shown to satisfy
\begin{equation}\label{CR_star_prod}T_f T_g \sim \sum_{j \ge 0} T_{\widehat C_j(f,g)},\end{equation}
where $\widehat C_j(f,g)$ is a classical pseudodifferential operator of order $-j$, $j \ge 0$. This means that for any $N\in \NN_0$, the operator $T_{f} T_g - \sum_{j=0}^N T_{\widehat C_j(f,g)}$ is of order $n-N-1$ (see Definition \ref{order_def}).
The pseudodifferential operators $\widehat C_j(f,g)$ are generally not determined uniquely, and their properties are less well understood than those of their analogues in the K\"{a}hler setting.

In \cite{gh2}, the authors define a notion of ``transversal Poisson bracket'', and show that the product of Toeplitz operators gives rise, through expansions such as (\ref{CR_star_prod}), to ``CR star products'' on certain suitable algebras of symbols. For instance, if $X$ admits a transversal CR $\RR$-action then it is possible to define a CR star product on $C^\infty(X)^{\RR}[[\hbar]]$, where $C^\infty(X)^{\RR}$ denotes the algebra of $\RR$-invariant functions. Notably, it is not assumed in \cite{gh2} that $X$ is strictly pseudoconvex (instead, a weaker non-degeneracy assumption is used).
The main result of the present paper (Theorem \ref{main_thm}) contributes towards the computation of the restriction of the principal symbol $c_{-1}(f,g)$ of $\widehat C_1(f,g)$ to the cone\footnote{Even after $\widehat C_0(f,g)$ is fixed, $c_{-1}(f,g)$ is only determined up to a homogeneous symbol $r$ of degree $-1$ such that $r|_{\Sigma} = 0$.}$$\Sigma = \{(x, -t \omega_0(x))\ | \ t > 0\}.$$
Fix $\widehat C_0(f,g) = \mcm_{fg}$ (the general case is similar). In a coordinate patch $(D,x)$, $$\begin{aligned} &T_{fg}(x,y) \equiv \int_0^\infty e^{it \phi(x,y)} b_{fg}(x,y,t) dt,\\ &T_{\widehat C_1(f,g)}(x,y) \equiv \int_0^\infty e^{it \phi(x,y)} b_{\widehat C_1(f,g)}(x,y,t) dt,\end{aligned}$$ where $\phi$ is the phase function specified in Theorem \ref{bdms_thm} and $$\begin{aligned} &b_{fg}(x,y,t) \sim \sum_{j \ge 0} b_{fg, j}(x,y) t^{n-j} \ \text{in } S^n_{1,0}(D \times D \times \RR_+),\\ &b_{\widehat C_1(f,g)}(x,y,t) \sim \sum_{j \ge 0} b_{\widehat C_1(f,g),j}(x,y) t^{n-j-1} \ \text{in } S^{n-1}_{1,0}(D \times D \times \RR_+).\end{aligned}$$Also, using largely the same arguments that are used in Sect. \ref{proof_sect}, $$\begin{aligned} &(T_f T_g)(x,y) \equiv \int_0^\infty e^{it \phi(x,y)} d_{f,g}(x,y,t) dt,\\ &d_{f,g}(x,y,t) \sim \sum_{j \ge 0} d_{f,g,j} t^{n-j}\ \text{in } S^n_{1,0}(D \times D \times \RR_+).\end{aligned}$$
The operator $T_f T_g - T_{fg}$ is of order $n-1$, and accordingly $$d_{f,g,0}(x,x) = b_{fg,0}(x,x) = \frac{f(x) g(x)}{2\pi^{n+1}}.$$
Now, since $T_f T_g - T_{fg} - T_{\widehat C_1(f,g)}$ is of order $n-2$, under suitable assumptions (see Lemma \ref{pointwise_uniqueness} or Corollary \ref{uniqueness_corollary}) it holds that $$\frac{c_{-1}(f,g)(x, -\omega_0(x))}{2\pi^{n+1}} =  b_{\widehat C_1(f,g),0}(x,x)= d_{f,g,1}(x,x) - b_{fg, 1}(x,x).$$
By homogeneity, this equation determines $c_{-1}(f,g)$ on the entire cone $\Sigma$. Thus, an explicit expression for $c_{-1}(f,g)\big{|}_{\Sigma}$ could be obtained through a computation of $d_{f,g,1}(x,x)$ and $b_{fg,1}(x,x)$.
At this point, we note that by Theorem \ref{main_thm}, $$b_{fg,1}(x,x) = \frac{R(x) f(x) g(x) - \square_b(fg)(x)}{4\pi^{n+1}}.$$
The computation of $d_{f,g,1}(x,x)$ is a natural continuation of the present work. In light of Proposition \ref{product_expansion}, it may be shown that $d_{f,g,1}(x,x)$ admits a formula involving, in particular, the subleading coefficients $b_{f,1}(x,x)$, $b_{g,1}(x,x)$ in the expansions associated with the distributional kernels of $T_f$, $T_g$. Consequently, Theorem \ref{main_thm} can also play a role in the computation of $d_{f,g,1}(x,x)$.
\section{Preliminaries}\label{prelim_sect}
\subsection{Notations}\label{notations}
Throughout this work, we use the notations specified in this section. We denote the set of non-negative integers by $\NN_0$. The set of positive real numbers will be denoted by $\RR_+$. If $\alpha = (\alpha_1, ..., \alpha_d)\in \NN_0^d$ is a multi-index, then
\begin{equation*}  \alpha! = \alpha_1! \cdots\alpha_d!,\ |\alpha| = \alpha_1 + \alpha_2 +\cdots+ \alpha_d.\end{equation*}
Let $x = (x_1, ..., x_d)$ be the standard coordinates on $\RR^d$. Then
\begin{equation*} \begin{aligned} &x^\alpha = x_1^{\alpha_1} \cdot ... \cdot x_d^{\alpha_d},\\
& \partial_{x_j} = \frac{\partial}{\partial x_j},\ j =1,...,d,\\ &\partial_x^\alpha = \partial^{\alpha_1}_{x_1} ... \partial^{\alpha_d}_{x_d} = \frac{\partial^{|\alpha|}}{\partial x^\alpha}.\end{aligned} \end{equation*}
In the specific case of second order derivatives, we also use the notation $$\partial_{x_j, x_k} = \partial_{x_j} \partial_{x_k} = \frac{\partial^2}{\partial x_j \partial x_k}.$$

If $U \subset \RR^d$ is open and $\phi : U \to \CC$ is differentiable, then we denote $$\phi'_x = (\partial_{x_1} \phi, ..., \partial_{x_d} \phi).$$%
If $\phi \in C^\infty(U)$, an almost analytic extension of $\phi$ is a function $ \widetilde \phi \in C^\infty(U_{\CC})$, where $U_{\CC} \subset \CC^d$ is open and $U_{\CC} \cap \RR^d = U$, such that $\widetilde \phi|_U = \phi$ and the functions $\partial_{\bar z_j} \widetilde \phi$, $j=1,...,d$, vanish to infinite order along $U$ (see \cite[Sect. 1]{ms}). Here, $z = x+iy \in U_\CC$ and $$\partial_{\bar z_j} = \frac 1 2 (\partial_{x_j} + i\partial_{y_j}).$$
For $f \in C^\infty(U \times U)$, we write $f(x,y) = \mco(|x-y|^\infty)$ if $$\partial^\alpha_x \partial^\beta_y f(x,x) = 0$$ for all $\alpha, \beta \in \NN_0^d$ and $x \in U$.

Let us recall some notions from microlocal analysis.
The space of smooth, compactly supported functions on an open set $U \subset \RR^d$ is denoted by $C^\infty_0(U)$, and $\mcd'(U)$ is the space of distributions on $U$. 
Let $A : C_0^\infty(V) \to \mcd'(U)$ be a continuous operator with Schwartz kernel $A(x,y) \in \mcd'(U \times V)$. We say that $A$ is a smoothing operator if $A(x,y) \in C^\infty(U \times V)$. Two continuous operators $A,B : C^\infty_0(V) \to \mcd'(U)$ are called equivalent if $A-B$ is a smoothing operator. We denote this equivalence by $$A(x,y) \equiv B(x,y) \mod C^\infty(U \times V),$$ or 
\[A\equiv B \mod C^\infty(U\times V)\]
or simply by $A(x,y) \equiv B(x,y)$ or $A\equiv B$.

The notion of H\"{o}rmander symbol spaces is central in this work.
\begin{definition} Let $D \subset \RR^{2n+1}$ be an open set and let $m \in \RR$. The H\"{o}rmander symbol space $S^m_{1,0}(D \times D \times \RR_+)$ consists of all functions $a \in C^\infty(D \times D \times \RR_+)$ which satisfy that for any compact set $K \Subset D \times D$, multi-indices $\alpha, \beta \in \NN_0^{2n+1}$ and $\gamma \in \NN_0$ there exists a constant $C_{K, \alpha, \beta, \gamma} > 0$ such that
\begin{equation*} \big{|}\partial^\alpha_x \partial^\beta_y \partial^\gamma_t a(x,y,t) \big{|} \le C_{K,\alpha, \beta, \gamma} (1+|t|)^{m - \gamma}\ \text{for all }(x,y,t) \in K \times \RR_+,\ t \ge 1.\end{equation*} Define also $S^{-\infty}_{1,0}(D \times D \times \RR_+) = \cap_{m \in \NN} S^{-m}_{1,0}(D \times D \times \RR_+)$.\end{definition}
If $(a_j)_{j \in \NN_0}$ is a sequence such that $a_j \in S^{m_j}_{1,0}(D \times D \times \RR_+)$, where $m_j \to -\infty$ as $j \to \infty$, then there exists $a \in S^{m_0'}_{1,0}(D \times D \times \RR_+)$, $m'_0 = \max_{j \ge 0} m_j$, such that 
\begin{equation*} a-\sum_{j=0}^{l-1} a_j \in S^{m_l'}_{1,0}(D \times D \times \RR_+),\ m'_l = \max_{j \ge l} m_j,\end{equation*}
for all $l=1,2,\ldots$. 
Moreover, $a$ is unique modulo $S^{-\infty}_{1,0}(D \times D \times \RR_+)$, and we will say that $$a \sim \sum^{+\infty}_{j = 0}a_j \ \text{in }S^{m_0'}_{1,0}(D \times D \times \RR_+).$$
Throughout, we will mostly deal with \textit{classical symbols}, defined as follows.
\begin{definition}\label{clhs} The space of classical symbols $S^m_{\cl}(D \times D \times \RR_+)$ consists of $a\in S^m_{1,0}(D \times D \times \RR_+)$ of the form
\begin{equation*} a(x,y,t) \sim \sum^{+\infty}_{l=0} a_l(x,y) t^{m-l} \text{in }S^m_{1,0}(D \times D \times \RR_+).\end{equation*}\end{definition}

The classical symbols appear in the representations of distribution kernels of Toeplitz operators, as well as more general types of operators.
\begin{definition}\label{order_def} A continuous operator $H: C^\infty(X) \to C^\infty(X)$ with distribution kernel $H(x,y)\in \mcd'(X \times X)$ is called a complex Fourier integral operator of Szeg\H{o} type of order $m \in \RR$ if it is smoothing away from the diagonal, and for any coordinate patch $(D,x)$ there exists $h \in S^m_{\cl}(D \times  D \times \RR_+)$ such that $$H(x,y) \equiv \int_0^\infty e^{it \phi(x,y)}h(x,y,t) dt,$$
where $\phi \in C^\infty(D \times D)$ is the phase function of Theorem \ref{bdms_thm}.
\end{definition}
We also consider pseudodifferential operators and their symbols.
\begin{definition} Let $D \subset \RR^{2n+1}$ be an open set, and let $m \in \RR$. The space of symbols $S^m_{1,0}(D\times \RR^{2n+1})$ consists of functions $p(x, \xi) \in C^\infty(D \times \RR^{2n+1})$ such that for any multi-indices $\alpha, \beta \in \NN_0^{2n+1}$ and compact set $K \subset D$ there exists a constant $C_{K, \alpha, \beta}$ such that
\begin{equation*} \big{|}\partial^\alpha_\xi \partial^\beta_x p(x, \xi) \big{|} \le C_{K, \alpha, \beta} (1+|\xi|)^{m-|\alpha|}\ \text{for all }x \in K,\ \xi \in \RR^{2n+1}.\end{equation*}
As before, we set $S^{-\infty}_{1,0}(D\times \RR^{2n+1}) = \cap_{m \in \NN} S^{-m}_{1,0}(D\times \RR^{2n+1})$. \end{definition}
The space of classical symbols $S^m_{\cl}(D\times \RR^{2n+1})$ consists of $p \in S^m_{1,0}(D\times \RR^{2n+1})$ admitting an expansion $$p(x,\xi) \sim \sum_{j \ge 0} p_j(x,\xi)$$ with $p_j$ positively homogeneous of degree $m-j$ (in the $\xi$ variable) when $|\xi| \ge 1$. 
We define $S^m_{1,0}(D\times D\times\mathbb R^{2n+1})$, $S^m_{{\rm cl\,}}(D\times D\times\mathbb R^{2n+1})$, $S^{-\infty}_{1,0}(D\times D\times\mathbb R^{2n+1})$ in a similar way. 
The space of classical pseudodifferential operators of order $m$ on $D$, denoted by $L^m_{\cl}(D)$, consists of pseudodifferential operators whose symbols belong to $S^m_{\cl}(D\times \RR^{2n+1})$. Similarly, if $X$ is a smooth manifold then $L^m_{\cl}(X)$ denotes the space of classical pseudodifferential operators on $X$.

We conclude this subsection by formulating a notion of subprincipal symbol of a classical pseudodifferential operator on a manifold equipped with a positive $s$-density, as follows.
\begin{lemma}\label{density_subprinc_symb} Let $X$ be a closed smooth manifold of dimension $n$, equipped with a positive $s$-density $\mu$, where $s \ne 0$, and let $E \in L^m_{\cl}(X)$. For each coordinate patch $(D, x)$ on $X$, write $$\mu\big{|}_D = \lambda |dx_1 \wedge ...\wedge dx_n|^s,\ \lambda \in C^\infty(D).$$
Let $e\in S^m_{\cl}(D\times \RR^n)$ be the total symbol of $E$ in the coordinate patch $(D,x)$, and write $$e(x,\xi) \sim \sum_{j \ge 0} e_j(x,\xi).$$
Define $e^{(D, x)}_{\sub}\in C^\infty(D \times \RR^n)$ by $$e^{(D,x)}_{\sub}(x,\xi) = e_1(x,\xi) + \frac i 2 \sum_{j=1}^n \partial_{x_j, \xi_j} e_0(x,\xi) +\frac i {2s}\sum_{j=1}^n \partial_{\xi_j} e_0(x,\xi) \partial_{x_j}(\log \lambda)(x).$$
Then the locally defined functions $e^{(D,x)}_{\sub}$ determine a well-defined global function $$e_{\sub} \in C^\infty(T^*X \setminus 0),$$
where $0 \subset T^*X$ denotes the zero section. We call $e_{\sub}$ the subprincipal symbol of $E$ with respect to $\mu$.
\end{lemma}
\begin{proof}
We identify $D$ with $x(D) \subset \mathbb R^n$, and we have that $$\mu\big{|}_D = \lambda|dx_1 \wedge ... \wedge dx_n|^s,\ \lambda\in C^\infty(D).$$ Let $\kappa : D \to D_\kappa$ be a diffeomorphism and put $y = \kappa(x)$. Let $$\mu_\kappa = (\kappa^{-1})^* \big{(}\mu\big{|}_D\big{)},$$
and write $$\mu_\kappa = \lambda_\kappa |dy_1 \wedge ... \wedge dy_n|^s,\ \lambda_\kappa \in C^\infty(D_\kappa).$$ Let $\Lambda(x) = \det(\kappa'(x))$.
Then the functions $\lambda$ and $\lambda_\kappa$ are related by $$\lambda_\kappa(\kappa(x)) = \frac {\lambda(x)} {|\Lambda(x)|^s}.$$ We consider $E$ as an operator $E \in L^m_{\cl}(D)$, and define $E_\kappa \in L^m_{\cl}(D_\kappa)$ by $$E_\kappa(u) = \left(E(u \circ \kappa) \right) \circ \kappa^{-1}.$$ Let $e_\kappa \in S^m_{\cl}(D_\kappa\times \RR^n)$ denote the total symbol of $E_\kappa$, and write $$e_\kappa \sim \sum_{j \ge 0} e_{\kappa, j}.$$
Define $$e_{ \sub}^{(D_\kappa, y)}(y, \eta) = e_{\kappa, 1}(y,\eta) + \frac i 2 \sum_{j=1}^n \partial_{y_j, \eta_j} e_{\kappa, 0}(y, \eta) + \frac i {2s}\sum_{j=1}^n \partial_{\eta_j} e_{\kappa,0}(y, \eta) \partial_{y_j}(\log \lambda_\kappa)(y).$$ In order to prove that $e_{\sub}$ is well-defined as a function on $T^*X\setminus 0$, it suffices to show that $$e^{(D_\kappa, y)}_{\sub}(\kappa(x), \eta) = e^{(D,x)}_{\sub}(x, (\kappa'(x))^t \eta),$$
where $(\kappa'(x))^t$ is the transpose of $\kappa'(x)$.

Denote $$\begin{aligned} &a(x,\xi) = \frac i {2s}\sum_{j=1}^n \partial_{\xi_j} e_0(x, \xi) \partial_{x_j} (\log \lambda)(x),\\ &a_\kappa(y, \eta) = \frac i {2s}\sum_{j=1}^n \partial_{\eta_j} e_{\kappa,0}(y, \eta) \partial_{y_j}(\log \lambda_\kappa)(y).\end{aligned}$$
Note that $$e_{\kappa, 0}(y, \eta) = e_0(x, (\kappa'(x))^t \eta).$$
Writing $\xi = (\xi_1, ..., \xi_n) = (\kappa'(x))^t \eta$, we have that (see \cite[p. 83]{hormander3}) $$\xi_l = \frac{\partial (y \cdot \eta) }{\partial x_l},\ \partial_{\eta_j} \xi_l =\partial_{x_l} y_j .$$
Thus,
$$\partial_{\eta_j} e_{\kappa,0} = \sum_{l=1}^n \partial_{\xi_l} e_0 \partial_{x_l} y_j.$$
Next, regarding $x = \kappa^{-1}(y)$ as a function of $y$, \begin{multline*}\partial_{y_j}(\log \lambda_{\kappa}) = \partial_{y_j}\big{(}\log \lambda - s \log(|\Lambda|)\big{)} = \frac{\partial_{y_j} \lambda }{\lambda} - s \frac{\partial_{y_j} (|\Lambda|)}{|\Lambda|} =\\ \sum_{m=1}^n \partial_{y_j} x_m \left[ \frac{\partial_{x_m} \lambda}{\lambda} - s \frac{\partial_{x_m}|\Lambda|}{|\Lambda|} \right] = \sum_{m=1}^n \partial_{y_j} x_m \left[\frac{\partial_{x_m} \lambda}{\lambda} - s \frac{\partial_{x_m} \Lambda}{\Lambda} \right].\end{multline*}
Thus,
\begin{multline*} \sum_{j=1}^n \partial_{\eta_j} e_{\kappa,0} \partial_{y_j}(\log\lambda_\kappa) = \\ \sum_{j=1}^n  \left(\sum_{l=1}^n \partial_{x_l} y_j \partial_{\xi_l} e_{0}\right) \left(\sum_{m=1}^n \partial_{y_j} x_m \left(\frac{\partial_{x_m} \lambda} \lambda - s \frac{\partial_{x_m}\Lambda}{\Lambda} \right) \right)\\=\sum_{l=1}^n \partial_{\xi_l} e_0 \sum_{j,m = 1}^n \partial_{x_l} y_j \partial_{y_j} x_m \left(\frac{\partial_{x_m} \lambda}{\lambda} - s\frac{\partial_{x_m} \Lambda}{\Lambda}\right) \\= \sum_{l=1}^n \partial_{\xi_l} e_0 \left(\frac{\partial_{x_l} \lambda}{\lambda} - s \frac{\partial_{x_l} \Lambda}{\Lambda} \right) = \sum_{l=1}^n \partial_{\xi_l} e_0 \left(\partial_{x_l}(\log \lambda) - s\frac{\partial_{x_l}\Lambda}{\Lambda}\right).\end{multline*}
Hence, $$ a_{\kappa}(\kappa(x), \eta) = a\left(x, (\kappa'(x))^t \eta\right)-\frac i 2 \sum_{l=1}^n\partial_{\xi_l} e_0\left(x, (\kappa'(x))^t \eta\right) \frac{\partial_{x_l} \Lambda(x)}{\Lambda(x)}.$$
In light of the transformation formula \cite[(18.1.33)]{hormander3}, it follows that $$e^{(D_\kappa, y)}_{\sub}(\kappa(x), \eta) = e^{(D, x)}_{\sub}(x, (\kappa'(x))^t \eta),$$
as required.
\end{proof}
\begin{remark} The subprincipal symbol defined in Lemma \ref{density_subprinc_symb} coincides with the usual subprincipal symbol (see \cite[Theorem 18.1.33]{hormander3}) of the operator $\mu^{\frac 1 {2s}} E \mu^{-\frac 1 {2s}}$, which acts on half-densities. \end{remark}
\subsection{Abstract CR manifolds and the Szeg\"{o} kernel}\label{CR_prelims}

Let $(X, T^{1,0}X)$ be a compact CR manifold of dimension $2n+1$, $n \ge 1$. This means that $T^{1,0}X \subset \CC TX$ is an involutive subbundle of rank $n$ such that for all $x\in X$, $$T^{1,0}_xX\cap T^{0,1}_xX = \{0\},$$ where $T^{0,1}X = \overline{T^{1,0}X}$. Assume also that $X$ is orientable. Then (\cite[1.1.2]{drago_tomas}) there exists a differential form $\omega_0 \in \Omega^1(X)$ satisfying
\begin{equation}\label{contact_form}\ker(\omega_0) = \re\big{(}T^{1,0}X \oplus T^{0,1}X\big{)}.\end{equation}
The \textit{Levi form} $\mcl_x : T_x^{1,0}X \times T_x^{1,0}X \to \CC$ with respect to $\omega_0$ can be specified by
\begin{equation}\label{levi_form} \mcl_x(u,v) = \frac i {2} d\omega_0(x) (u, \bar v),\ u,v \in T_x^{1,0}X.\end{equation}
Throughout, we assume that $\mcl_x$ is positive definite for all $x \in X$, in which case $X$ is called \textit{strongly pseudoconvex} and $\omega_0$ is a contact form. The Reeb vector field associated with $\omega_0$ is the real vector field $\mct \in C^\infty(X, TX)$, defined by the equations
\begin{equation*} \omega_0(\mct) = -1,\ \iota_{\mct} d\omega_0 = 0.\end{equation*}

The Levi form gives rise to a Hermitian metric $\langle \cdot, \cdot \rangle_{\mcl}$ on $\CC TX$ called the \textit{Levi metric}, which is specified by
\begin{equation}\label{levi_metric} \begin{aligned} &\langle u, v \rangle_{\mcl} = \mcl_x(u, v)\ \text{for } u,v \in T^{1,0}_x X,\\ &\langle \bar u, \bar v \rangle_{\mcl} = \overline{\langle u, v \rangle}_{\mcl} \ \text{for } u,v \in T^{1,0}_x X,\\ & T^{1,0}_xX \perp T^{0,1}_xX,\ \mct \perp (T^{1,0}X \oplus T^{0,1}X),\ \langle \mct, \mct \rangle_{\mcl} = 1.\end{aligned} \end{equation} 
Let $T^{*1,0}X$ and $T^{*0,1}X$ denote the subbundles of $\CC T^*X$ consisting of forms which annihilate $\CC \mct \oplus T^{0,1}X$ and $\CC \mct \oplus T^{1,0}X$, respectively.
The Levi metric yields (by duality) a Hermitian metric on $\CC T^*X$ such that
\begin{equation*} \CC T^*X = T^{*1,0}X \oplus T^{*0,1}X \oplus \CC \omega_0\end{equation*}
is an orthogonal direct sum decomposition.

Let $\Omega^{0,1}(X):= C^\infty(X, T^{*0,1}X).$
The \textit{tangential Cauchy-Riemann operator} \begin{equation}\label{tcro}  \bar \partial_b : C^\infty(X) \to \Omega^{0,1}(X)\end{equation} is defined as
\begin{equation*} \bar \partial_b = \pi^{(0,1)} \circ d,\end{equation*}
where $d : C^\infty(X) \to \Omega^1(X)$ is the exterior derivative and $$\pi^{(0,1)} : \CC T^*X \to T^{*0,1}X$$ is the fiberwise orthogonal projection.

Let $dv_X$ be the volume form associated with the Levi metric. The Hermitian metric on $\CC T^*X$ gives rise to an inner product $(\,\cdot\,|\,\cdot\,)$ on $\Omega^{0,1}(X)$ (by integration with respect to $dv_X$). Let $L^2_{(0,1)}(X)=L^2_{(0,1)}(X, dv_X)$ be the completion of $\Omega^{0,1}(X)$ with respect to $(\,\cdot\,|\,\cdot\,)$. We also write $(\,\cdot\,|\,\cdot\,)$ to denote the inner product on $C^\infty(X)$ induced by the volume form $dv_X$. Let $L^2(X)=L^2(X, dv_X)$  be the completion of $C^\infty(X)$ with respect to $(\,\cdot\,|\,\cdot\,)$. We extend $\bar\partial_b$ to $L^2$ space:
\[\bar\partial_b: {\rm Dom\,}\bar\partial_b\subset L^2(X)\rightarrow L^2_{(0,1)}(X),\]
where ${\rm Dom\,}\bar\partial_b:=\{u\in L^2(X);\, \bar\partial_bu\in L^2_{(0,1)}(X)\}$. We also write 
\[\bar\partial^*_b: {\rm Dom\,}\bar\partial^*_b\subset L^2_{(0,1)}(X)\rightarrow L^2(X)\]
to denote the $L^2$ adjoint of $\bar\partial_b$.

The \textit{Kohn Laplacian} $\square_b : C^\infty(X) \to C^\infty(X)$ is defined by $$\square_b = \bar \partial_b^* \bar \partial_b,$$ where $\bar \partial_b^*$ is the formal adjoint of $\bar \partial_b$. We extend $\square_b$ to $L^2$ space:

\[\square_b: {\rm Dom\,}\square_b\subset L^2(X)\rightarrow L^2(X),\]
where ${\rm Dom\,}\square_b:=\{u\in L^2(X);\, u\in{\rm Dom\,}\bar\partial_b, \bar\partial_bu\in{\rm Dom\,}\bar\partial^*_b\}$.

Throughout this work, we assume that the range of $\square_b$ is closed.

The \textit{Szeg\H{o} projection} is the orthogonal projection $$\Pi : L^2(X, dv_X) \to \ker(\bar \partial_b).$$ The distribution kernel $$\Pi(x,y) \in \mcd'(X \times X)$$ of the Szeg\H{o} projection is known as the \textit{Szeg\H{o} kernel}. The microlocal behavior of $\Pi(x,y)$ is addressed in Theorems \ref{bdms_thm}, \ref{hsiao_shen_thm_intro}.

Finally, the main result of the present paper is formulated using the so-called \textit{Tanaka-Webster scalar curvature}, which may be specified as follows.
\begin{definition}[\cite{tanaka}]\label{tan_web_con} The Tanaka-Webster connection is the unique affine connection $\nabla$ such that
\begin{enumerate}
\item{$\nabla_U V \in C^\infty(X, HX)$ for every $U \in C^\infty(X, TX)$ and $V \in C^\infty(X, HX)$, where $HX$ is the Levi distribution (\ref{levi_distribution}).}
\item{$\nabla \mct = \nabla J = \nabla d\omega_0 = 0$, where $J$ is the complex structure (\ref{cplx_str}) on $HX$.}
\item{The torsion $\tau$ of $\nabla$ satisfies $$\tau(U,V) = -d\omega_0(U,V) \mct,\ \tau(\mct, J U) = -J \tau(\mct, U),\ U, V \in C^\infty(X, HX).$$}
\end{enumerate}
\end{definition}
The Tanaka-Webster connection is compatible with the Levi metric (\ref{levi_metric}). Next, let $\{Z_j\}_{j=1}^n$ be a local frame of $T^{1,0}X$ with dual frame $\{\theta^j\}_{j=1}^n$, and denote $$Z_{\bar j} = \overline{Z_j},\ \theta^{\bar j} = \overline{\theta^j},\ j=1,...,n.$$ Then $$\nabla Z_j = \sum_{k=1}^n \omega_j^k \otimes Z_k,\ \nabla Z_{\bar j} = \sum_{k = 1}^n \omega_{\bar j}^{\bar k} \otimes Z_{\bar k},$$
and $(\omega_j^k)_{j,k=1,...,n}$ is the connection one form with respect to the chosen local frame. The Tanaka-Webster curvature two form $(\Theta_j^k)_{j, k = 1,...,n}$ is defined by $$\Theta_j^k = d\omega_j^k - \sum_{l = 1}^n \omega_j^l \wedge \omega_l^k.$$
A direct computation implies that $$\Theta_j^k = C_{0,j}^k \wedge \omega_0 + \sum_{l,m=1}^n\left[ R^k_{jl \bar m} \theta^l \wedge \theta^{\bar m} + A^k_{jlm} \theta^l \wedge \theta^m + B^k_{jlm} \theta^{\bar l} \wedge \theta^{\bar m}\right],$$
where $C_{0,j}^k$ is some differential one-form. The terms $R^k_{j l \bar m}$ are the components of the \textit{pseudohermitian curvature tensor} associated with $\nabla$, and $$R_{j \bar m} = \sum_{k = 1}^n R^k_{j k \bar m}$$
is called the \textit{pseudohermitian Ricci curvature}. Finally, write $$-d\omega_0 = i \sum_{j,k} g_{j \bar k} \theta^j \wedge \theta^{\bar k},$$
and let $g^{\bar l m}$ denote the components of the inverse of the matrix $(g_{j \bar k})_{j,k=1,...,n}$.
\begin{definition} The \textit{Tanaka-Webster scalar curvature} $R \in C^\infty(X)$ with respect to $\omega_0$ is specified by \begin{equation}\label{tw_curv}R = \sum_{j, m = 1}^n g^{\bar m j} R_{j \bar m}.\end{equation} \end{definition}
\subsection{Local representations}
The computations underlying the main results of this work are facilitated by a convenient choice of local coordinates. Adapting the construction in the proof of \cite[Lemma 3.2]{hormanderpaper} and the computations of \cite[Sect.~3.2]{hsiao_shen}, we obtain the following.
\begin{proposition}\label{canonic_coor} Assume that the range of $\square_b$ is closed. Fix $p \in X$. There exist local coordinates $x= (x_1, x_2, ..., x_{2n+1})$ on an open neighborhood $D$ of $p$, with $x(p) = 0$, which satisfy the following properties.
\begin{enumerate}
\item{Denote $z_j = x_{2j-1} + i x_{2j}$. Then $\omega_0 \big{|}_D$ is specified by
\begin{equation*} \omega_0(x) = (1+q_2(x))d x_{2n+1} + \frac i 2 \sum_{j=1}^n \left(\bar z_j dz_j - z_j d\bar z_j\right) + \mco(|x|^3),\end{equation*}
where $q_2$ is a quadratic form satisfying $\partial_{z_j} \partial_{\bar z_k} q_2 = 0$ for all $j,k = 1,...,n$. Here $\partial_{z_j} = \frac 1 2 \big{(}\partial_{x_{2j-1}} - i \partial_{x_{2j}}\big{)}$ and $\partial_{\bar z_k} = \frac 1 2\big{(}\partial_{x_{2k-1}}+i \partial_{x_{2k}}\big{)}$. }
\item{There exists a local frame $Z_1, ..., Z_n$ of $T^{1,0}X\big{|}_D$, orthonormal with respect to the Levi metric, such that
\begin{equation*} Z_j(x) = \sqrt 2 \left(\partial_{z_j} - \frac i 2 \bar z_j \partial_{x_{2n+1}} \right) + \mco(|x|^2).\end{equation*}}
\item{The vector field $\mct\big{|}_D$ satisfies $$\mct(x) = -\partial_{x_{2n+1}}+\frac 1 2 \sum_{j=1}^n \Big{(}\partial_{x_{2j}} q_2 \partial_{x_{2j-1}} - \partial_{x_{2j-1}} q_2 \partial_{x_{2j}} \Big{)} +\mco(|x|^2).$$}
\item{The volume form $dv_X\big{|}_D$ is given by
\begin{equation*} dv_X(x) = \lambda(x) dx_1 \wedge ... \wedge dx_{2n+1},\end{equation*}
with
\begin{equation*} \lambda(0) = 1,\ \partial_{x_j} \lambda(0) = 0,\ j=1,...,2n+1.\end{equation*}}
\end{enumerate}
\end{proposition}
The phase function and symbol in Theorem \ref{bdms_thm} are not determined uniquely. We will make use of the following choices.
\begin{theorem}[\cite{hsiao_shen}]\label{hsiao_shen_thm} Let $(D,x)$ be the coordinates of Proposition \ref{canonic_coor}. Then in Theorem \ref{bdms_thm}, we may choose $$\phi \in C^\infty(D \times D),\ a \in S^n_{\cl}(D \times D \times \RR_+)$$ such that $$a_0(x,y) = \frac 1 {2\pi^{n+1}} + \mco(|x-y|^\infty),$$ and $\phi$ satisfies (\ref{phi_props}), and additionally, writing $w_j = y_{2j-1} + i y_{2j}$,%
\begin{equation}\label{phi_a_forms}\phi(x,y) = -x_{2n+1} + y_{2n+1} + \frac i 2 \sum_{j=1}^n \left(|z_j - w_j|^2 + (\bar z_j w_j - z_j \bar w_j) \right) + \mco(|(x,y)|^3).\end{equation}
\end{theorem}

Finally, we note that the Kohn Laplacian also has a convenient expression in terms of our chosen coordinates.
\begin{lemma} In the coordinates of Proposition \ref{canonic_coor}, the Kohn Laplacian satisfies for all $f\in C^\infty(X)$, 
\begin{equation*} \square_bf(0) = -\frac 1 2\sum_{j=1}^{2n} \partial_{x_j}^2 f(0) -in \partial_{x_{2n+1}}f(0).\end{equation*} \end{lemma}
\begin{proof}
A straightforward computation yields (cf. \cite[pp. 154-156]{beals_greiner})
\begin{equation*} \square_bf(0) = \sum_{j=1}^n \overline{Z_j}^* \overline{Z_j} f(0),\end{equation*}
and since $\text{div}\left(\partial_{\bar z_j} +\frac i 2  z_j \partial_{x_{2n+1}}\right)(0)=0$, so that $\overline{Z_j}^* = - Z_j + \mco(|x|)$, we obtain the required. \end{proof}

\subsection{Some geometry on $T^*X$}\label{diffr_op_def}
In what follows, we specify a differential operator $P: C^\infty(T^*X) \to C^\infty(T^*X)$ such that for any $F \in C^\infty(T^*X)$, in the local coordinates of Proposition \ref{canonic_coor},
$$P(F)(0, -\omega_0(0)) = \sum_{j=1}^n \left[\big{(}\partial_{x_{2j}, \xi_{2j-1}} F\big{)}(0, -\omega_0(0)) - \big{(}\partial_{x_{2j-1}, \xi_{2j}} F \big{)}(0, -\omega_0(0))\right].$$

Let $\omega=\omega_{T^*X}$ denote the standard symplectic form on $T^*X$. Let $(D, x)$ be a coordinate patch as in Proposition \ref{canonic_coor}, and let $(x, \xi)$ denote the induced coordinates on $T^*D$. Then $$\omega_{T^*X} = \sum_{j=1}^{2n+1} dx_j \wedge d\xi_j.$$
The Hamiltonian vector field of $F \in C^\infty(T^*D)$ is the unique vector field $X_F$ such that $$\omega_{T^*X}(X_F, \cdot) = -dF.$$
A straightforward computation produces $$X_F = \sum_{j=1}^{2n+1} \left( \partial_{x_j} F \partial_{\xi_j}-\partial_{\xi_j} F \partial_{x_j}\right).$$
Let $\Div : \Vect(T^*X) \to C^\infty(T^*X)$ denote the divergence with respect to the Liouville volume form $\frac{\omega^{\wedge (2n+1)}}{(2n+1)!} $, where $\Vect(T^*X):=C^\infty(T^*X,TT^*X)$. Since $\omega_{T^*X}\big{|}_{T^*D}$ is the standard symplectic form, the divergence is specified by $$\Div \left(\sum_{j=1}^{2n+1}\left[ X_j \partial_{x_j} + \Xi_j \partial_{\xi_j} \right]\right) = \sum_{j=1}^{2n+1}\left( \partial_{x_j} X_j + \partial_{\xi_j} \Xi_j\right).$$

Let $HX$ denote the Levi distribution of $X$, that is, $HX \subset TX$ is the unique subbundle such that $$\CC HX = T^{1,0}X \oplus T^{0,1}X.$$
Then $HX$ carries a natural complex structure $J : HX \to HX$ given by \begin{equation}\label{cplx_str}J(u + \bar u) = iu - i \bar u,\ u \in T^{1,0}X.\end{equation}
Let $(D,x)$ be a coordinate patch as in Proposition \ref{canonic_coor}. Introduce a local frame $X_j$, $j = 1,...,2n$, of $HX$ by\begin{equation}\label{levi_frame}\begin{aligned} & X_{2j-1} = \frac 1{\sqrt 2} \left(Z_j + \bar Z_j\right) = \partial_{x_{2j-1}} - x_{2j} \partial_{x_{2n+1}} + R_{2j-1},\\ & X_{2j} = \frac 1{\sqrt 2} \left(i Z_j + \overline{i Z_j} \right) = \partial_{x_{2j}} + x_{2j-1} \partial_{x_{2n+1}} + R_{2j},\\ &R_j = \mco(|x|^2),\ j = 1,...,2n.\end{aligned}\end{equation}
Then \begin{equation*} J X_{2j-1} = X_{2j},\ J X_{2j} = - X_{2j-1}.\end{equation*}
Additionally, write \begin{equation}\label{x2n+1}X_{2n+1} = - \mct = \partial_{x_{2n+1}} + \sum_{j=1}^{2n}u_j(x) \partial_{x_j}+ R_{2n+1},\end{equation} where $R_{2n+1} = \mco(|x|^2)$ and \begin{equation}\label{u_j_formula}u_j = \left\{\begin{array}{ll} \frac 1 2 \partial_{x_{2m-1}} q_2 & \text{if }j = 2m,\\[3pt] -\frac 1 2 \partial_{x_{2m}} q_2 & \text{if }j = 2m-1,\end{array}\right.\end{equation}
where $q_2(x)$ is the quadratic form specified in Proposition \ref{canonic_coor}.

\begin{lemma}\label{christoffel_symbols_lemma} Consider the Christoffel symbols $\{\Gamma_{jk}^l\}_{j,k,l = 1,...,2n+1}$ specified by $$\nabla_{\partial_{x_j}} d{x_k} = \sum_{l=1}^{2n+1} \Gamma^l_{jk} d x_l.$$
Then $$\Gamma^l_{jk} = \left\{\begin{array}{ll} -1+ \mco(|x|)& \text{if } k = 2n+1,\ j = 2m,\ l = 2m-1,\\[2pt] 1+ \mco(|x|)& \text{if } k = 2n+1,\ j= 2m-1,\ l = 2m,\\[2pt] \partial_{x_j} u_k + \mco(|x|) & \text{if } k\le 2n,\ l = 2n+1,\\[2pt]  \mco(|x|) & \text{otherwise} \end{array}\right. ,$$ 
where $\nabla$ is the Tanaka-Webster connection given by Definition~\ref{tan_web_con}. 
\end{lemma}
\begin{proof}
First, we note that for $j, k = 1,...,n$, \begin{equation}\label{zero_der}\nabla_{\bar Z_k} Z_j(0) = \nabla_{Z_k} Z_j(0) = \nabla_{\mct} Z_j(0) = 0.\end{equation}
Indeed, since $\nabla$ preserves $HX$ and $\nabla J = 0$, it holds that $\nabla$ preserves $T^{1,0}X$, namely, $\nabla_Z Z_j \in T^{1,0}X$ for any $Z \in C^\infty(X, \CC TX)$. Thus, using the properties of the torsion $\tau$ of $\nabla$, we see that $$\nabla_{\bar Z_k} Z_j = \big{(}[\bar Z_k, Z_j] \big{)}^{1,0},$$ where the latter is the component of $[\bar Z_k, Z_j]$ lying in $T^{1,0}X$. Now, a quick computation gives $$\big{(}[\bar Z_k, Z_j]\big{)}^{1,0}(0) = 0,$$ as required. Next, $\nabla$ is compatible with the Levi metric (and $Z_j,\ j = 1,...,n$ is an orthonormal frame), hence at $p = 0$, $$0 = Z_k \langle Z_j, Z_l\rangle_{\mcl} = \langle \nabla_{Z_k} Z_j, Z_l \rangle_{\mcl} + \langle Z_j, \nabla_{\bar Z_k} Z_l \rangle_{\mcl} = \langle \nabla_{Z_k} Z_j, Z_l \rangle_{\mcl},$$
and since $\nabla_{Z_k} Z_j \in C^\infty(X, T^{1,0}X)$, we conclude that $\nabla_{Z_k} Z_j(0) = 0$. Similarly, it is readily verified that $\tau(\mct, Z_j) \in T^{0,1}X$, which implies that $\nabla_{\mct} Z_j = \big{(} [\mct, Z_j] \big{)}^{1,0}$, and a direct computation produces $\nabla_{\mct} Z_j(0) = 0$.

Now, noting that equation (\ref{zero_der}) implies that $$\nabla_{\bar Z_k} \bar Z_j(0) = \nabla_{Z_k} \bar Z_j(0) = \nabla_{\mct} \bar Z_j(0) = 0,\ j,k = 1,...,n,$$
we conclude that $$\nabla_{X_k} X_j(0) = 0,\ j,k = 1,...,2n+1.$$

Next, we are interested in $$\nabla_{\partial_{x_j}} \partial_{x_k} = \sum_{l = 1}^{2n+1} \tilde \Gamma^l_{jk} \partial_{x_l}.$$
For $j = 1,...,2n$, set $$\nu_j = \left\{\begin{array}{ll} x_{2k}& \text{if } j = 2k-1,\\ -x_{2k-1} & \text{if }j = 2k.\end{array}\right.$$ Then using (\ref{levi_frame}), $$\partial_{x_j} = X_j + \nu_j X_{2n+1}  - R_j.$$
Thus, for $j,k = 1,...,2n$,
\begin{multline*} \nabla_{\partial_{x_j}} \partial_{x_k} = \nabla_{X_j + \nu_j X_{2n+1} - R_j} (X_k + \nu_k X_{2n+1} - R_k) = \\ \nabla_{X_j}(X_k + \nu_k X_{2n+1} -  R_k)+ \nu_j \nabla_{X_{2n+1}}(X_k + \nu_k X_{2n+1} - R_k)- \nabla_{ R_j}(X_k + \nu_k X_{2n+1} -  R_k)\\ = \nabla_{X_j}(X_k + \nu_k X_{2n+1} - R_k) + \mco(|x|).\end{multline*}
Next,
\begin{multline*} \nabla_{X_j}(X_k + \nu_k X_{2n+1} - R_k) = \nabla_{X_j} X_k  + \nu_k \nabla_{X_j} X_{2n+1} + X_j(\nu_k) X_{2n+1} - \nabla_{X_j} R_k\\ = X_j(\nu_k) X_{2n+1} + \mco(|x|) = \left(\partial_{x_j} \nu_k \right)\partial_{x_{2n+1}} + \mco(|x|).\end{multline*}
Noting that $$\partial_{x_m} \nu_k = \left\{\begin{array}{ll} 1 & \text{if } m = 2j,\ k = 2j-1,\\ -1 & \text{if } m = 2j-1,\ k = 2j,\\ 0 & \text{otherwise} \end{array}\right.,$$
we conclude that for $m, k = 1,...,2n$, $$\nabla_{\partial_{x_m}}\partial_{x_k} = \left\{\begin{array}{ll}\partial_{x_{2n+1}} + \mco(|x|) & \text{if } m = 2j,\ k = 2j-1,\\ -\partial_{x_{2n+1}} + \mco(|x|) & \text{if } m = 2j-1,\ k = 2j,\\ \mco(|x|)& \text{otherwise} \end{array}\right. .$$
Similarly, noting that $$X_{2n+1} = \partial_{x_{2n+1}} + \sum_{k=1}^{2n} u_k(x) \partial_{x_k} + R_{2n+1},$$ we obtain for $j = 1,...,2n+1$, $$\nabla_{\partial_{x_j}} \partial_{x_{2n+1}} = -\sum_{k=1}^{2n} \partial_{x_j} u_k(x) \partial_{x_k} + \mco(|x|).$$
In a similar way we obtain that $\nabla_{\partial_{x_{2n+1}}} \partial_{x_j} = \mco(|x|)$, $j = 1,...,2n$.
Finally, the required follows from the fact that since $dx_j$, $j = 1,...,2n+1$ is the dual frame of $\partial_{x_j}$, $j = 1,...,2n+1$, hence $\Gamma^l_{jk} = - \tilde \Gamma^k_{jl}$.
\end{proof}
Let $\pi_{T^*X} : T^*X \to X$ be the canonical projection and let $$\Ver(T^*X) = \ker d\pi_{T^*X}\subset T T^*X$$ denote the vertical subbundle of $TT^*X$.
The Tanaka-Webster connection gives rise to a direct sum decomposition $$T T^*X = \Hor(T^*X) \oplus \Ver(T^*X),$$
where $\Hor(T^*X)$ is the horizontal subbundle of $T^*X$. Then $\Hor(T^*X)$ may be identified with the pullback bundle $\pi_{T^*X}^* TX$, allowing us to lift the Levi distribution $HX \subset TX$ to obtain a distribution $H X_{\Hor} \subset \Hor(T^*X)$. Let $J_{\Hor} : HX_{\Hor}\to HX_{\Hor}$ be the lift of the complex structure $J : HX \to HX$. In what follows, the horizontal lift of a vector field $V$ on $X$ is denoted by $V^{\Hor}$, and $u_k(x)$, $k = 1,...,n$, are as specified in (\ref{u_j_formula}).
\begin{corollary}\label{hor_frame} A standard computation produces for $j = 1,...,n$,$$\begin{aligned} &\partial_{x_{2j}}^{\Hor} = \partial_{x_{2j}} + \xi_{2n+1} \partial_{\xi_{2j-1}}+ K_{2j}(\xi)\partial_{\xi_{2n+1}} + r_{2j},\\ &\partial_{x_{2j-1}}^{\Hor} = \partial_{x_{2j-1}} -\xi_{2n+1} \partial_{\xi_{2j}} + K_{2j-1}(\xi) \partial_{\xi_{2n+1}} + r_{2j-1},\\ &\partial_{x_{2n+1}}^{\Hor} = \partial_{x_{2n+1}} + K_{2n+1}( \xi) \partial_{\xi_{2n+1}} + r_{2n+1},\end{aligned}$$
where for $m = 1,...,2n+1$, it holds that $r_m = \mco(|x|)$, $r_m \in \Ver(T^*D)$ and $$K_m(-\omega_0(0)) =0,\ \partial_{\xi_{2n+1}} K_m = 0.$$ Note that in the right hand side above, we view $\partial_{x_j}$ as a vector field on $T^*D$. Consequently, for $j = 1,...,n$, $$\begin{aligned} &X_{2j}^{\Hor} = \partial^{\Hor}_{x_{2j}} + x_{2j-1} \partial_{x_{2n+1}}^{\Hor} + \mco(|x|^2),\\ &X_{2j-1}^{\Hor} = \partial_{x_{2j-1}}^{\Hor} - x_{2j} \partial_{x_{2n+1}}^{\Hor}  + \mco(|x|^2),\end{aligned}$$
and additionally $$X_{2n+1}^{\Hor} = \partial_{x_{2n+1}}^{\Hor}  +  \sum_{k=1}^{2n} u_k(x) \partial_{x_k}^{\Hor}+ \mco(|x|^2).$$ \end{corollary}
\begin{proof}
Recall that the horizontal lift of a vector $V\in T_x X$ to $T_{(x,\xi)} T^*X$ is defined as follows. First, consider any curve $\gamma : (-\varepsilon, \varepsilon) \to X$ such that $$\gamma(0) = x,\ \dot \gamma(0) = V.$$ By parallel transport (defined through the Tanaka-Webster connection), $\gamma$ lifts to a curve $\gamma_{T*X} : (-\varepsilon, \varepsilon) \to T^*X$ such that $\gamma_{T^*X}(0) = (x, \xi)$. The horizontal lift $V^{\Hor}$ of $V$ is then given by $V^{\Hor} = \dot \gamma_{T^*X}(0)$.
In local coordinates, it may be verified that the horizontal lift of $$V = \sum_{j=1}^{2n+1} a_j \partial_{x_j}$$ is given by \begin{equation}\label{horlift_formula} V^{\Hor} = \sum_{j=1}^{2n+1} a_j \partial_{x_j}^{\Hor},\end{equation}
where \begin{equation}\label{horlift_coor_vec_fields}\partial_{x_j}^{\Hor} = \partial_{x_j} - \sum_{l,k = 1}^{2n+1} \xi_k \Gamma^l_{jk} \partial_{\xi_l}.\end{equation}
Substituting the formulas for the Christoffel symbols specified in Lemma \ref{christoffel_symbols_lemma} into equation (\ref{horlift_coor_vec_fields}), we obtain $$\begin{aligned} &\partial_{x_{2j}}^{\Hor} = \partial_{x_{2j}} - \sum_{l,k = 1}^{2n+1} \xi_k \Gamma^l_{2j, k} \partial_{\xi_l} = \partial_{x_{2j}} + \xi_{2n+1} \partial_{\xi_{2j-1}}+K_{2j}(\xi)\partial_{\xi_{2n+1}} + r_{2j},\\ & K_{2j}(\xi) =-\sum_{k=1}^{2n} \xi_k \partial_{x_{2j}} u_k =  \frac 1 2 \sum_{k=1}^n \big{(}\xi_{2k-1} \partial_{x_{2j}, x_{2k}} q_2 - \xi_{2k} \partial_{x_{2j}, x_{2k-1}} q_2\big{)},\end{aligned}$$
and similarly $$\begin{aligned} &\partial_{x_{2j-1}}^{\Hor} = \partial_{x_{2j-1}} - \xi_{2n+1} \partial_{\xi_{2j}} + K_{2j-1}(\xi)\partial_{\xi_{2n+1}} + r_{2j-1},\\ &\partial_{x_{2n+1}}^{\Hor} = \partial_{x_{2n+1}} + K_{2n+1}(\xi)\partial_{\xi_{2n+1}}+ r_{2n+1},\end{aligned}$$
where $$r_m \in \Ver(T^*D),\ r_m = \mco(|x|),\ m = 1,...,2n+1,$$ and $$K_m(\xi) =  -\sum_{k=1}^{2n} \xi_{k} \partial_{x_m} u_k.$$
Finally, applying equation (\ref{horlift_formula}) to the vector fields $X_1,...,X_{2n+1}$ (noting their expressions (\ref{levi_frame}), (\ref{x2n+1}) in the local frame $\partial_{x_1},...,\partial_{x_{2n+1}}$), we obtain the required formulas for $X_1^{\Hor},...,X_{2n+1}^{\Hor}$.
\end{proof}
\begin{lemma}\label{hor_levi_comp} Let $V$ be a vector field on $T^*X$ such that in the coordinates of Proposition \ref{canonic_coor},$$V = \sum_{j=1}^{2n+1} (a_j \partial_{x_j} + b_j \partial_{\xi_j}).$$
Let $V^{HX_{\Hor}}$ denote the component of $V$ which lies inside $HX_{\Hor}$. Then $$V^{H X_{\Hor}} = \sum_{j=1}^{2n} \big{(}a_j- a_{2n+1} u_j\big{)} X_j^{\Hor} + \mco(|x|^2),$$
with $u_j$, $j =1,...,2n$, specified in (\ref{u_j_formula}).\end{lemma}
\begin{proof}
Let $V^{\Hor}$ and $V^{\Ver}$ denote the horizontal and vertical components of $V$. Write $$V = \sum_{j=1}^{2n+1} a_j \partial_{x_j}^{\Hor} + \sum_{j=1}^{2n+1}a_j (\partial_{x_j} - \partial_{x_j}^{\Hor}) + \sum_{j=1}^{2n+1} b_j \partial_{\xi_j}.$$
Note that $\partial_{x_j} - \partial_{x_j}^{\Hor}$ is a vertical vector field, which means that $$V^{\Hor} =\sum_{j=1}^{2n+1} a_j \partial_{x_j}^{\Hor}.$$
Now, since $u_j(x) = \mco(|x|)$, we obtain\begin{multline*} V^{\Hor} = \sum_{j=1}^{2n+1} a_j X_j^{\Hor} + \sum_{j=1}^n \big{(} a_{2j}(\partial_{x_{2j}}^{\Hor} - X_{2j}^{\Hor}) + a_{2j-1}(\partial_{x_{2j-1}}^{\Hor} - X_{2j-1}^{\Hor})\big{)} \\+ a_{2n+1}(\partial_{x_{2n+1}}^{\Hor} - X_{2n+1}^{\Hor})=\sum_{j=1}^{2n+1} a_j X_j^{\Hor} \\+ \Bigg{[}\sum_{j=1}^n\big{(}a_{2j-1} x_{2j} - a_{2j} x_{2j-1}\big{)}\Bigg{]} X_{2n+1}^{\Hor} - a_{2n+1} \sum_{j=1}^{2n} u_j X_j^{\Hor}+\mco(|x|^2)\\ = \sum_{j=1}^{2n} \big{(}a_j - a_{2n+1} u_j\big{)}X_j^{\Hor} + \Bigg{[} a_{2n+1} + \sum_{j=1}^n \big{(} a_{2j-1} x_{2j} - a_{2j} x_{2j-1}\big{)}\Bigg{]} X_{2n+1}^{\Hor} + \mco(|x|^2).\end{multline*} 
\end{proof}
\begin{lemma}Let $V$ be a vector field on $T^*X$ such that in the coordinates of Proposition \ref{canonic_coor},$$V = \sum_{j=1}^{2n+1}( a_j \partial_{x_j} + b_j \partial_{\xi_j}).$$ Then \begin{multline*}\Div\left(J_{\Hor} V^{HX_{\Hor}}\right)(0, -\omega_0(0)) =\\ \sum_{j=1}^n \big{[}\partial_{x_{2j}} a_{2j-1}(0, -\omega_0(0)) - \partial_{x_{2j-1}} a_{2j}(0,-\omega_0(0))\big{]} - \sum_{j=1}^{2n} \partial_{\xi_j} a_j(0, -\omega_0(0)).\end{multline*}\end{lemma}
\begin{proof}
By Lemma \ref{hor_levi_comp}, $$V^{HX_{\Hor}} = \sum_{j=1}^{2n} (a_j - a_{2n+1} u_j) X_j^{\Hor} + \mco(|x|^2).$$ Thus, $$J_{\Hor} V^{HX_{\Hor}} = \sum_{j=1}^n \big{(}(a_{2j-1}- a_{2n+1} u_{2j-1}) X_{2j}^{\Hor} - (a_{2j}-a_{2n+1} u_{2j}) X_{2j-1}^{\Hor}\big{)} + \mco(|x|^2).$$
For a vector field of the form $a X_{2j-1}^{\Hor}$, using the divergence product formula, $$\Div(a X_{2j-1}^{\Hor}) = X_{2j-1}^{\Hor}(a) + a \Div(X_{2j-1}^{\Hor}).$$
Using Corollary \ref{hor_frame}, \begin{multline*}X_{2j-1}^{\Hor} = \partial_{x_{2j-1}} - x_{2j} \partial_{x_{2n+1}} - \xi_{2n+1} \partial_{\xi_{2j}} +(K_{2j-1}(\xi) - x_{2j} K_{2n+1}(\xi)) \partial_{\xi_{2n+1}} \\ + r_{2j-1} + \mco(|x|^2),\end{multline*}
where $r_{2j-1} = \mco(|x|)$ is a vertical term. Thus, noting that $-\omega_0(0) = (0, ..., 0, -1)$, we immediately obtain that $$\Div(X_{2j-1}^{\Hor})(0, -\omega_0(0)) =0,$$
and additionally \begin{equation}\label{odd_div}X_{2j-1}^{\Hor}(a)(0, -\omega_0(0)) = \partial_{x_{2j-1}}a(0, -\omega_0(0)) + \partial_{\xi_{2j}} a(0,-\omega_0(0)).\end{equation}
Similarly, for a vector field of the form $a X_{2j}^{\Hor}$, \begin{equation}\label{even_div}\Div(a X_{2j}^{\Hor})(0, -\omega_0(0)) = \partial_{x_{2j}} a(0, -\omega_0(0)) - \partial_{\xi_{2j-1}} a(0, -\omega_0(0)).\end{equation}
Thus, using (\ref{odd_div}), (\ref{even_div}), we obtain the desired formula by a straightforward computation of $$\Div\Bigg{[} \sum_{j=1}^n \big{(}a_{2j-1} X_{2j}^{\Hor} - a_{2j} X_{2j-1}^{\Hor}\big{)}\Bigg{]}(0, -\omega_0(0)).$$
Finally, denoting $$\hat U = \sum_{j=1}^n u_{2j} X_{2j-1}^{\Hor} - u_{2j-1} X_{2j}^{\Hor},$$
it holds that (by the properties of $q_2$ specified in Proposition \ref{canonic_coor}) $$\Div \hat U(0, -\omega_0(0))= \sum_{j=1}^n \big{(}\partial_{x_{2j-1}} u_{2j}(0) -\partial_{x_{2j}} u_{2j-1}(0) \big{)} =   2 \sum_{j=1}^n \partial_{z_j, \bar z_j} q_2(0) = 0,$$
and since $u_j = \mco(|x|)$, $j=1,...,2n$, we conclude (by the product formula) that $\Div(a_{2n+1} \hat U)(0, -\omega_0(0)) = 0$. 
\end{proof}
Let $\omega_1, ..., \omega_{2n}$ be the frame of $H^* X = \re(T^{*1,0}X \oplus T^{*0,1}X)$ which is dual to $X_1,..., X_{2n}$. Then $$\omega_j = dx_j -u_j dx_{2n+1} + \mco(|x|^2),\ j = 1,...,2n,$$
where $u_j$ is specified in (\ref{u_j_formula}), $u_j(x) = \mco(|x|)$. By a slight abuse of notation, we let $J : H^* X \to H^*X$ denote the dual of $J : HX \to HX$. Then $$J \omega_{2j-1} = -\omega_{2j},\ J \omega_{2j} = \omega_{2j-1}.$$
The vertical bundle $\Ver(T^*X)\subset T T^*X$ is naturally isomorphic to $\pi^*_{T^*X} T^*X$. Let $\omega^{\Ver} \in \Ver(T^*X)$ denote the vertical lift of $\omega \in T^*X$. Note that $$dx_j^{\Ver} = \partial_{\xi_j},\ j = 1,...,2n+1.$$
Let $H X_{\Ver} \subset \Ver(T^*X)$ be the vertical lift of $H^*X$, and $J_{\Ver} : HX_{\Ver} \to HX_{\Ver}$ the lift of $J : H^*X \to H^*X$. 
\begin{lemma} Let $V$ be a vector field on $T^*X$ such that in the coordinates of Proposition \ref{canonic_coor},$$V = \sum_{j=1}^{2n+1} (a_j \partial_{x_j} + b_j \partial_{\xi_j}).$$ Consider the decomposition $$T T^*X = HX_{\Hor} \oplus \RR \mct^{\Hor} \oplus HX_{\Ver} \oplus \RR \omega_0^{\Ver},$$
and let $V^{HX_{\Ver}}$ denote the component of $V$ lying in $HX_{\Ver}$. Then \begin{multline*} V^{HX_{\Ver}} = \sum_{j=1}^n \big{[}(b_{2j} + \xi_{2n+1} a_{2j-1}) \omega^{\Ver}_{2j} +(b_{2j-1} - \xi_{2n+1} a_{2j}) \omega_{2j-1}^{\Ver} \big{]} + \mco(|x|).\end{multline*} \end{lemma}
\begin{proof}
Writing $$V = \sum_{j=1}^{2n+1} a_j \partial_{x_j}^{\Hor} + \sum_{j=1}^{2n+1}a_j (\partial_{x_j} - \partial_{x_j}^{\Hor}) + \sum_{j=1}^{2n+1} b_j \partial_{\xi_j},$$ it is clear that with respect to the decomposition $T T^*X = \Hor(T^*X) \oplus \Ver(T^*X)$, the vertical component of $V$ is given by  $$\hat V = \sum_{j=1}^{2n+1} b_j \partial_{\xi_j} - \sum_{j=1}^{2n+1} a_j(\partial_{x_j}^{\Hor} - \partial_{x_j}).$$
That is, \begin{multline*} \hat V = \sum_{j=1}^{2n+1} b_j \partial_{\xi_j} - \sum_{j=1}^n \big{(} a_{2j}(\partial_{x_{2j}}^{\Hor} - \partial_{x_{2j}}) + a_{2j-1}(\partial_{x_{2j-1}}^{\Hor} - \partial_{x_{2j-1}})\big{)}\\ - a_{2n+1}(\partial_{x_{2n+1}}^{\Hor} - \partial_{x_{2n+1}}) = \sum_{j=1}^{2n+1} b_j \partial_{\xi_j}-\\  \sum_{j=1}^n \big{(} a_{2j}(\xi_{2n+1} \partial_{\xi_{2j-1}} + K_{2j} \partial_{\xi_{2n+1}} + r_{2j}) + a_{2j-1}(-\xi_{2n+1} \partial_{\xi_{2j}} + K_{2j-1} \partial_{\xi_{2n+1}})\big{)} \\- a_{2n+1}K_{2n+1}\partial_{\xi_{2n+1}} - r_{2n+1} \\ = \sum_{j=1}^{2n+1} b_j \partial_{\xi_j} + \xi_{2n+1} \sum_{j=1}^n \big{(} a_{2j-1} \partial_{\xi_{2j}} - a_{2j} \partial_{\xi_{2j-1}}\big{)} - \Big{[} \sum_{j=1}^{2n+1} a_j K_j\Big{]} \partial_{\xi_{2n+1}} + r .\end{multline*}
Here, $ r$ is vertical and satisfies $ r = \mco(|x|)$.
Noting that $$\omega_j^{\Ver} = \partial_{\xi_j} + \hat r_j,\ j = 1,...,2n,$$ with $\hat r_j$ vertical and $\hat r_j = \mco(|x|)$, and similarly $\omega_0^{\Ver} = \partial_{\xi_{2n+1}} + \hat r_0$, $\hat r_0$ vertical and $\hat r_0 = \mco(|x|)$, we obtain that $$\hat V = \sum_{j=1}^{2n} b_j \omega_j^{\Ver} + \xi_{2n+1} \sum_{j=1}^n \big{(} a_{2j-1} \omega_{2j}^{\Ver} - a_{2j} \omega_{2j-1}^{\Ver} \big{)} +\Big{[}b_{2n+1} - \sum_{j=1}^{2n+1} a_j K_j \Big{]} \omega_0^{\Ver} + \hat r,$$
where $\hat r$ is vertical, $\hat r = \mco(|x|)$. From this, the required readily follows.
\end{proof}
\begin{lemma}  Let $V$ be a vector field on $T^*X$ such that in the coordinates of Proposition \ref{canonic_coor},$$V = \sum_{j=1}^{2n+1}( a_j \partial_{x_j} + b_j \partial_{\xi_j}).$$ Consider the decomposition $$T T^*X = HX_{\Hor} \oplus \RR \mct^{\Hor} \oplus HX_{\Ver} \oplus \RR \omega_0^{\Ver},$$
and let $V^{HX_{\Ver}}$ denote the component of $V$ lying in $HX_{\Ver}$. Then \begin{multline*} \Div(J_{\Ver}V^{HX_{\Ver}})(0, -\omega_0(0)) = \\ \sum_{j=1}^n \big{[} \partial_{\xi_{2j-1}} b_{2j}(0, -\omega_0(0)) - \partial_{\xi_{2j}} b_{2j-1}(0, -\omega_0(0)) \big{]} -\sum_{j=1}^{2n} \partial_{\xi_j} a_j(0, -\omega_0(0)).\end{multline*} \end{lemma}
\begin{proof}
In light of the previous lemma,
$$V^{HX_{\Ver}} = \sum_{j=1}^n \big{[}(b_{2j} + \xi_{2n+1} a_{2j-1}) \omega^{\Ver}_{2j} +(b_{2j-1} - \xi_{2n+1} a_{2j}) \omega_{2j-1}^{\Ver} \big{]} + \mco(|x|),$$
hence
$$J_{\Ver}V^{HX_{\Ver}} = \sum_{j=1}^n \big{[}(b_{2j} + \xi_{2n+1} a_{2j-1}) \omega^{\Ver}_{2j-1} +(\xi_{2n+1} a_{2j}-b_{2j-1}) \omega_{2j}^{\Ver} \big{]} + \mco(|x|),$$
that is $$J_{\Ver} V^{HX_{\Ver}} =\sum_{j=1}^n \big{[}(b_{2j} + \xi_{2n+1} a_{2j-1}) \partial_{\xi_{2j-1}} +(\xi_{2n+1} a_{2j}-b_{2j-1}) \partial_{\xi_{2j}} \big{]} + \mco(|x|).$$
Thus, \begin{multline*} \Div(J_{\Ver}V^{HX_{\Ver}}) = \\\sum_{j=1}^n \big{[} \partial_{\xi_{2j-1}} b_{2j} + \xi_{2n+1} \partial_{\xi_{2j-1}} a_{2j-1} + \xi_{2n+1} \partial_{\xi_{2j}} a_{2j} - \partial_{\xi_{2j}} b_{2j-1} \big{]} + \mco(|x|),\end{multline*}
and finally
\begin{multline*} \Div(J_{\Ver}V^{HX_{\Ver}})(0, -\omega_0(0)) = \\\sum_{j=1}^n \big{[} \partial_{\xi_{2j-1}} b_{2j}(0, -\omega_0(0)) - \partial_{\xi_{2j}} b_{2j-1}(0, -\omega_0(0)) \big{]} -\sum_{j=1}^{2n} \partial_{\xi_j} a_j(0, -\omega_0(0)).\end{multline*}
\end{proof}
Finally, define a morphism $J_{T^*X} : T T^*X \to T T^*X$ by \begin{equation}\label{lifted_J}\begin{aligned} &J_{T^*X} \big{|}_{HX_{\Hor}} = J_{\Hor},\\ &J_{T^*X} \big{|}_{HX_{\Ver}} = - J_{\Ver},\\ &J_{T^*X} \mct^{\Hor} = J_{T^*X} \omega_0^{\Ver} = 0. \end{aligned}\end{equation}
\begin{corollary}\label{div_jhamvf} Define a differential operator $P : C^\infty(T^*X) \to C^\infty(T^*X)$ by $$P(F) = -\frac 1 2\Div(J_{T^*X} X_F).$$
Then in the coordinates of Proposition \ref{canonic_coor}, $$P(F)(0, -\omega_0(0)) = \sum_{j=1}^n\big{[}\partial_{x_{2j}, \xi_{2j-1}} F(0,-\omega_0(0)) - \partial_{x_{2j-1}, \xi_{2j}} F(0, -\omega_0(0))\big{]}.$$
\end{corollary}
\begin{proof}
Let $V$ be a vector field on $T^*X$ such that
$$V = \sum_{j=1}^{2n+1} \left(a_j \partial_{x_j} + b_j \partial_{\xi_j}\right).$$
Then (using the notations above) $$J_{T^*X} V = J_{\Hor} V^{HX_{\Hor}} - J_{\Ver} V^{HX_{\Ver}},$$
and consequently $$\begin{aligned} &\Div(J_{T^*X} V)(0, -\omega_0(0)) =\\ &\sum_{j=1}^n \big{[}\partial_{x_{2j}} a_{2j-1}(0, -\omega_0(0)) - \partial_{x_{2j-1}} a_{2j}(0,-\omega_0(0))\big{]} - \sum_{j=1}^{2n} \partial_{\xi_j} a_j(0, -\omega_0(0))\\ &-\sum_{j=1}^n \big{[} \partial_{\xi_{2j-1}} b_{2j}(0, -\omega_0(0)) - \partial_{\xi_{2j}} b_{2j-1}(0, -\omega_0(0)) \big{]} +\sum_{j=1}^{2n} \partial_{\xi_j} a_j(0, -\omega_0(0))\\ &=\sum_{j=1}^n \big{[}\partial_{x_{2j}} a_{2j-1}(0, -\omega_0(0)) - \partial_{x_{2j-1}} a_{2j}(0,-\omega_0(0))\big{]}\\ &-\sum_{j=1}^n \big{[} \partial_{\xi_{2j-1}} b_{2j}(0, -\omega_0(0)) - \partial_{\xi_{2j}} b_{2j-1}(0, -\omega_0(0)) \big{]}.\end{aligned}$$
Inserting $$a_j = -\partial_{\xi_j} F,\ b_j = \partial_{x_j} F,$$
we obtain the required.
\end{proof}
\section{The expansion of the kernel of a Toeplitz operator}\label{proof_sect}
Let $E\in L^m_{\cl}(X)$ denote a classical pseudodifferential operator of order $m \in \mathbb R$.
Let $T_E(x,y) \in \mcd'(X \times X)$ denote the Schwartz kernel of the Toeplitz operator $T_E = \Pi E \Pi$. Then $T_E(x,y)$ is smooth away from the diagonal (by standard results on wavefront sets of products). Let $Q_E(x,y)$ denote the Schwartz kernel of the operator $Q_E = E \Pi$. In what follows, we study $T_E(x,y)$ by using the relation
\begin{equation}\label{kern_comp} T_E(x,y) \equiv \int_X \Pi(x,u) Q_E(u, y) dv_X(u).\end{equation}
The leading term in the expansion of the symbol of $\Pi(x,y)$ has been computed in the seminal paper \cite{bdms}. The computations in \cite{bdms} can also be used in order to derive the leading term in the expansion of the symbol of $T_E(x,y)$, as follows (cf. \cite[Theorem 4.4]{gh2}).
\begin{lemma}\label{lead_factor} Let $(D, x)$ be an open coordinate patch such that $$\Pi(x,y) \equiv \int_0^\infty e^{it \phi(x,y)} a(x,y,t) dt,$$
where $\phi$, $a$ are as in Theorem \ref{bdms_thm}. Let $E\in L^m_{\cl}(X)$ denote a classical pseudodifferential operator of order $m \in \mathbb R$ with principal symbol $e_0 \in C^\infty(T^*X \setminus 0)$. Then there exists $b \in S^{n+m}_{\cl}(D \times D \times \RR_+)$,
$$b \sim \sum_{j \ge 0} b_j t^{n+m-j} \ \text{in } S^{n+m}_{1,0}(D \times D \times \RR_+),$$
such that $$T_E(x,y) \equiv \int_0^\infty e^{it \phi(x,y)} b(x,y,t) dt,$$
and for any $x \in D$, $$b_0(x,x) = \frac{e_0(x, -\omega_0(x))}{2\pi^{n+1}}.$$
Finally, $\mct_y b_0 = 0$ on $D \times D$. \end{lemma}
\begin{proof}
The same arguments as in the proof of Proposition \ref{Q_E_prop} imply that $$Q_E(x,y) \equiv \int_0^\infty e^{it \phi(x,y)} c(x,y,t) dt,$$
where $c \in S^{n+m}_{\cl}(D \times D \times \RR_+)$ satisfies that $$c_0(x,x)= \frac{e_0(x, \phi'_x(x,x))}{2\pi^{n+1}} = \frac{e_0(x, -\omega_0(x))}{2\pi^{n+1}}.$$
Next, since $T_E = \Pi Q_E$, it follows from the stationary phase formula (\cite[Theorem 2.3]{ms}, cf. \cite[Lemma 5.3]{hsiaomari2}) that there exists a complex valued phase function $\varphi \in C^\infty(D \times D)$ satisfying (\ref{phi_props}) and $d \in S^{n+m}_{\cl}(D \times D \times \RR_+)$, $$ \begin{aligned} &d(x,y,t) \sim \sum_{j = 0}^\infty d_j(x,y) t^{n+m-j}\ \text{in } S^{n+m}_{1,0}(D \times D \times \RR_+),\\ & d_0(x,x) = 2 \pi^{n+1} a_0(x,x) c_0(x,x) = c_0(x,x), \end{aligned}$$such that $$T_E(x,y) \equiv \int_0^\infty e^{it \varphi(x,y)} d(x,y,t) dt.$$
Moreover, as shown in \cite[Theorem 5.4]{hsiaomari2} (see also \cite[Sect. 8]{hsiaomari2}), there exists $f \in C^\infty(D \times D)$ with $f(x,x) = 1$ such that $$\varphi(x,y) = f(x,y) \phi(x,y) + \mco(|x-y|^\infty).$$
Hence, we may assume without loss of generality that $\varphi(x,y) = f(x,y) \phi(x,y)$, and using Lemma \ref{phase_factorization}, conclude that $$T_E(x,y) \equiv \int_0^\infty e^{it f(x,y) \phi(x,y)} d(x,y,t) \equiv \int_0^\infty e^{it \phi(x,y)} b(x,y,t) dt,$$
with $b \in S^{n+m}_{\cl}(D \times D \times \RR_+)$ satisfying $$\begin{aligned} &b(x,y,t) \sim \sum_{j = 0}^\infty b_j(x,y) t^{n+m-j}\ \text{in } S^{n+m}_{1,0}(D \times D \times \RR_+),\\ &b_0(x,x) = \frac{c_0(x,x)}{(f(x,x))^{n+m+1}} = c_0(x,x).\end{aligned}$$
Finally, by applying the same arguments of \cite[Lemma 3.3]{hsiao_shen}, we may assume without loss of generality that $$\mct_y b_0(x,y) = 0$$ for every $(x,y) \in D \times D$.
\end{proof}
\begin{corollary} Let $(D, x)$ be an open coordinate patch. Let $\hat \phi \in C^\infty(D \times D)$ satisfy (\ref{phi_props}). Assume that $\hat \phi$ is equivalent in the sense of \cite{ms} to the phase function $\phi$ of Theorem \ref{bdms_thm}. Let $E\in L^m_{\cl}(X)$ be a classical pseudodifferential operator of order $m \in \mathbb R$ with principal symbol $e_0 \in C^\infty(T^*X \setminus 0)$. Then there exists $\hat b \in S^{n+m}_{\cl}(D \times D \times \RR_+)$,
$$\hat b \sim \sum_{j \ge 0} \hat b_j t^{n+m-j} \ \text{in } S^{n+m}_{1,0}(D \times D \times \RR_+),$$
such that $$T_E(x,y) \equiv \int_0^\infty e^{it \hat \phi(x,y)} \hat b(x,y,t) dt,$$
and additionally, for any $x \in D$, $$\hat b_0(x,x) = \frac{e_0(x, -\omega_0(x))}{2\pi^{n+1}}.$$
Finally, $\mct_y \hat b_0 = 0$ on $D \times D$. \end{corollary}
\begin{proof}
First, note that there exists $B \in S^{n+m}_{\cl}(D \times D \times \RR_+)$, $$B \sim \sum_{j \ge 0} B_j t^{n+m-j}\ \text{in } S^{n+m}_{1,0}(D \times D \times \RR_+),$$
such that for any $x \in D$, $$B_0(x,x) = \frac{e_0(x, -\omega_0(x))}{2\pi^{n+1}},$$
and $$T_E(x,y) \equiv \int_0^\infty e^{it \hat \phi(x,y)} B(x,y,t) dt.$$
Indeed, since $\hat \phi$ satisfies (\ref{phi_props}) and it is equivalent to $\phi$, there exists (see \cite[Sect. 8]{hsiaomari2}) $f \in C^\infty(D \times D)$ with $f(x,x) = 1$ such that $$\hat \phi(x,y) = f(x,y) \phi(x,y) + \mco(|x-y|^\infty).$$ Hence, we may assume without loss of generality that $\hat \phi(x,y) = f(x,y) \phi(x,y)$. Next, in light of Lemma \ref{lead_factor}, we conclude using Lemma \ref{phase_factorization} that there exists $B \in S^{n+m}_{\cl}(D \times D \times \RR_+)$ as required. Finally, using the same argument as in \cite[Lemma 3.3]{hsiao_shen}, it may be shown that we can replace $B \in S^{n+m}_{\cl}(D \times D \times \RR_+)$ with $\hat b \in S^{n+m}_{\cl}(D \times D \times \RR_+)$ as required.
\end{proof}
\subsection{Point-wise computations}
In this subsection, we fix an arbitrary point $p\in X$ and express $T_E(x,y)$ (locally) using a phase function and a symbol which are especially adapted so as to simplify computations at $p$. Then, we consider the first two coefficients in the expansion of the symbol, and evaluate them at $p$. In light of Lemma \ref{pointwise_uniqueness}, this will suffice to complete the proof of Theorem \ref{main_thm}. More precisely, fix $p \in X$ and let $(D, x)$ be a coordinate patch as in Proposition \ref{canonic_coor}, with $x(p) = 0$. Throughout, we identify $D$ with $x(D)$. Let $$\phi \in C^\infty(D \times D),\ a \in S^n_{\cl}(D \times D \times \RR_+)$$ be as in Theorem \ref{hsiao_shen_thm}. Then by Lemma \ref{lead_factor}, $$T_E(x,y) \equiv \int_0^\infty e^{it \phi(x,y)} b(x,y,t)dt,$$
where $$\begin{aligned} &b\in S^{n+m}_{\cl}(D \times D \times \RR_+),\\ &b \sim \sum_{j \ge 0} b_j t^{n+m-j}\ \text{in }S^{n+m}_{1,0}(D \times D \times \RR_+),\\ &\mct_y b_0 = 0\ \text{on }D \times D,\end{aligned}$$ and for all $x \in D$, $$b_0(x,x) = \frac{e_0(x, -\omega_0(x))}{2\pi^{n+1}}.$$
Since $\mct_y^2 \phi(0,0) = 0$, Lemma \ref{pointwise_uniqueness} implies that in order to complete the proof of Theorem \ref{main_thm}, it suffices to verify that $b_1(0,0)$ agrees with the expression specified in formula (\ref{main_formula}). In fact, our computation will not involve $\phi$ or $b$ directly. Rather, by Malgrange's preparation theorem (\cite[Theorem 7.5.5]{hormander1}), there exist $f, g \in C^\infty(D \times D)$ such that $f(x,y) = 1 + \mco(|(x,y)|^2)$ and \begin{equation}\label{phiPhi}\phi(x,y) = f(x,y)(y_{2n+1} + g(x,y')),\ y' = (y_1,...,y_{2n}).\end{equation}
Instead of $\phi$, we will prefer to represent $T_E(x,y)$ using the phase function \begin{equation}\label{loc_indep_phase}\Phi(x,y) = y_{2n+1} + g(x,y'),\end{equation} because of its simple dependence on $y_{2n+1}$. Note that $\mct_y^2 \Phi(0,0) = 0$. Hence, a slight modification of Lemma \ref{pointwise_uniqueness} implies that if $$\begin{aligned} &B \in S^{n+m}_{\cl}(D \times D \times \RR_+),\\ & B \sim \sum_{j \ge 0} B_j t^{n+m-j}\ \text{in } S^{n+m}_{1,0}(D \times D \times \RR_+) \end{aligned}$$ satisfies
$$\begin{aligned} & \int_0^\infty e^{it \phi(x,y)} b(x,y,t) dt \equiv \int_0^\infty e^{it \Phi(x,y)} B(x,y,t) dt\\ & B_0(0,0) = b_0(0,0),\\ &\mct_y B_0(0,0) = 0, \end{aligned}$$
Then $$B_1(0,0) = b_1(0,0).$$ Thus, the following result, which is the main result of the present subsection, suffices to complete the proof of Theorem \ref{main_thm}.
\begin{theorem}\label{main_thm_exp_sec} Fix $p \in X$. Let $(D,x)$ be a coordinate patch as specified in Proposition \ref{canonic_coor}, with $x(p) = 0$, and let $(x,\xi)$ be the induced coordinates on $T^*D$.  Let $E\in L^m_{\cl}(X)$ denote a classical pseudodifferential operator of order $m \in \mathbb R$ with principal symbol $e_0 \in C^\infty(T^*D \setminus 0)$, and define $\mce_0(x) = e_0(x, -\omega_0(x))$.  Let $e_{\sub} \in C^\infty(T^*X\setminus 0)$ be the subprincipal symbol of $E$ with respect to the $1$-density\footnote{Here $dv_X$ is the volume form specified in Sect. \ref{CR_prelims}.} $|dv_X|$, as specified in Lemma \ref{density_subprinc_symb}. Let $$P : C^\infty(T^*X) \to C^\infty(T^*X)$$ be the differential operator specified in Corollary \ref{div_jhamvf}. Let $\Phi \in C^\infty(D \times D)$ be the phase function specified in Proposition \ref{hsiao_shen_indep}. Then
\begin{equation*} T_E(x,y) \equiv \int_0^\infty e^{it  \Phi(x,y)} B(x,y,t) dt,\end{equation*}
where $B(x,y,t) \in S^{n+m}_{\cl}(D \times D \times \RR_+)$ is independent of $y_{2n+1}$ and satisfies
\begin{equation*} B(x,y,t) \sim \sum_{j=0}^\infty B_{j}(x,y) t^{n+m-j}\text{ in } S^{n+m}_{1,0}(D \times D \times \RR_+)\end{equation*}
with $B_{l} \in C^\infty(D \times D)$ independent of $y_{2n+1}$ for all $l \ge 0$, \begin{equation*}B_{0}(0,0) = \frac{e_0(0, -\omega_0(0))}{2\pi^{n+1}},\end{equation*}
and \begin{multline*} B_{1}(0,0) = \frac 1 {4\pi^{n+1}}\left[R(0) \mathcal E_0(0) -\square_b \mathcal E_0(0)  + P(e_0)(0, -\omega_0(0))\right]\\+\frac 1 {4\pi^{n+1}} \left[2 e_{\sub}(0, -\omega_0(0))-i m \mct \mce_0(0)\right]. \end{multline*}\end{theorem}
The proof of Theorem \ref{main_thm_exp_sec} relies (in particular) on the following properties of $\Phi$, which were mostly established in \cite{hsiao_shen} (see pages 352-354).
\begin{proposition}[\cite{hsiao_shen}]\label{hsiao_shen_indep} Let $(D,x)$ be a coordinate patch as in Proposition \ref{canonic_coor}. Let $\Phi \in C^\infty(D \times D)$ be the phase function as in (\ref{loc_indep_phase}). Then $$\Phi(x,y) = y_{2n+1}-x_{2n+1} + \frac i 2 \sum_{j=1}^n \left(|z_j - w_j|^2 + (\bar z_j w_j - z_j \bar w_j) \right) + \mco(|(x,y')|^4),$$
and \begin{equation}\label{generalized_phi_props}\begin{aligned} & \im(\Phi) \ge 0,\\ & \Phi(x,y) = 0 \ \text{if and only if }x=y,\\ & d_x\Phi(x,x) = - d_y \Phi(x,x) = -\nu(x)\omega_0(x)\end{aligned}\end{equation}
for a positive $\nu \in C^\infty(D)$ which satisfies $\nu(x) = 1 -q_2(x)+ \mco(|x|^3)$, where $q_2$ is as specified in Proposition \ref{canonic_coor}. Additionally, $$\Pi(x,y) \equiv \int_0^\infty e^{it \Phi(x,y)} A(x,y,t)dt,$$
where $$\begin{aligned} &A \in S^n_{\cl}(D \times D \times \RR_+),\\ &A\sim \sum_{j \ge 0} A_j t^{n-j}\ \text{in } S^n_{1,0}(D \times D \times \RR_+),\\ &\partial_{y_{2n+1}} A = \partial_{y_{2n+1}} A_j = 0\ \text{for all }j \ge 0,\end{aligned}$$ with $$\begin{aligned} &A_0(x,y') = \frac {1} {2\pi^{n+1}} + \mco(|(x,y')|^2),\\ &\partial_{z_j, \bar z_k} A_0(0,0) = \partial_{w_j, \bar w_k} A_0(0,0) = 0,\ j,k=1,...,n,\end{aligned}$$
and $$A_1(0,0) = \frac{R(0)}{4\pi^{n+1}}.$$ Here, $R$ denotes the Tanaka-Webster scalar curvature (\ref{tw_curv}) associated with $\omega_0$.\end{proposition}
\begin{proof}
Write $\phi(x,y) = \phi_0(x,y) + p_3(x,y) + \mco(|(x,y)|^4)$ with $$\phi_0(x,y) = -x_{2n+1} + y_{2n+1} + \frac i 2 \sum_{j=1}^n \left(|z_j - w_j|^2 + (\bar z_j w_j - z_j \bar w_j) \right)$$ and $p_3$ a homogeneous cubic term. For $j=1,..,n$, since $\omega_0(Z_j) = 0$, we obtain that $$Z_j = \sqrt 2(\partial_{z_j} -\frac i 2 \bar z_j \partial_{x_{2n+1}}) + r_j,$$
where the $\partial_{x_{2n+1}}$ coefficient of $r_j$ is of order $\mco(|x|^3)$. Hence, $$\overline{Z}_{j,x} \phi = \sqrt 2 \partial_{\bar z_j} p_3 + O(|(x,y)|^3),\ Z_{j,y} \phi = \sqrt 2 \partial_{w_j} p_3 + \mco(|(x,y)|^3),$$
and noting that $\overline{\partial}_{b,x}\phi$ and $\partial_{b,y}\phi$ vanish to infinite order on the diagonal, we conclude that $\partial_{\bar z_j} p_3 = \partial_{w_j} p_3 = 0$. Next, we note that $p_3(z, x_{2n+1}, \bar z, x_{2n+1}) = 0$ (since $\phi(x,x) = 0$), hence $p_3(z, x_{2n+1}, \bar w, y_{2n+1}) = 0$ whenever $x_{2n+1} = y_{2n+1}$, that is, $$p_3(x,y) = p_2(z, x_{2n+1},\bar w, y_{2n+1})(y_{2n+1} - x_{2n+1})$$ for a homogeneous quadratic polynomial $p_2$.  Thus, the equation $$\phi(x,y) = f(x,y) \Phi(x,y)$$ implies that $\Phi = \phi_0 + \mco(|(x, y')|^4)$, and $f(x,y) = 1 + p_2(x,y) + \mco(|(x,y)|^3)$. Also, noting that $\partial_{x_{2n+1}} \phi(x,x) = -1 - q_2(x)+ \mco(|x|^3)$, we conclude that $$\begin{aligned} &f(x,x) = 1 + q_2(x) + \mco(|x|^3),\\ &\partial_{z_j} \partial_{\bar z_k} f(0,0) = \partial_{w_j} \partial_{\bar w_k} f(0,0) = 0,\ j,k = 1,...,n.\end{aligned}$$
This implies that $A_0$ can be chosen such that it satisfies the properties above (see \cite{hsiao_shen}, pp. 352-354). Finally, the formula for $A_1(0,0)$ was obtained in \cite{hsiao_shen}, and it is also derived in Lemma \ref{computation_lemma} of the present paper (relying on computations from \cite{hsiao_shen}).
\end{proof}

Let $e(x,\xi)$ denote the total symbol of $E$ over $D$. Then
\begin{equation*} e(x,\xi) \sim \sum_{j \ge 0} e_j(x,\xi)\ \text{in } S^m_{\cl}(D\times \RR^{2n+1}),\end{equation*}
where $e_j$ is positively homogeneous of degree $m-j$ in the $\xi$ variable when $|\xi| \ge 1$. In what follows, we fix $\Phi$, $A$ as specified in Proposition \ref{hsiao_shen_indep}, such that $$\Pi(x,y) \equiv \int_0^\infty e^{it \Phi(x,y)}A(x,y,t) dt.$$
\begin{proposition}\label{Q_E_prop} There exists $C(x,y,t) \in S^{n+m}_{\cl}( D \times D \times \RR_+)$ which does not depend on $y_{2n+1}$ such that
\begin{equation*} \begin{aligned} &Q_E(x,y) \equiv \int_0^\infty e^{it \Phi(x,y)} C(x,y,t) dt,\\ &C(x,y,t) \sim \sum_{j =0}^\infty C_{j}(x,y) t^{n+m-j}\ \text{in } S^{n+m}_{1,0}(D \times D \times \RR_+),\end{aligned}\end{equation*}
and $C_{j}(x,y)$ are independent of $y_{2n+1}$ for all $j \in \NN_0$, and
\begin{equation*} \begin{aligned}&C_{0}(x,y) = \tilde e_0(x, \Phi_x'(x,y))A_0(x,y),\\ & \begin{aligned} C_{1}(x,y) = &\ \tilde e_0(x, \Phi'_x(x,y)) A_1(x,y) \\ &+\frac 1 {2\pi^{n+1}}\Big{(} \tilde e_1(x,  \Phi'_x(x,y)) - i \sum_{|\alpha| = 2} \frac{\partial^\alpha_{\xi} \tilde e_0(x, \Phi'_x(x,y)) \partial^\alpha_{x} \Phi(x,y)}{\alpha!}\Big{)}\\ &+\mco(|(x,y)|).\end{aligned}\end{aligned}\end{equation*}
Here $\tilde e_0$, $\tilde e_1$ are almost analytic extensions of $e_0$, $e_1$.
\end{proposition}
\begin{proof}
Note that since $\Pi(x,y)$ is smooth away from the diagonal, we can assume that $A$, $A_j$, $j \ge 0$, are properly supported on $D \times D$.
The existence of $$C(x,y,t) \in S^{n+m}_{\cl}(D \times D \times \RR_+)$$ such that $$Q_E(x,y) \equiv \int_0^\infty e^{i t \Phi(x,y)} C(x,y,t) dt$$is essentially explained in \cite[p. 178]{ms} (see \cite[Theorem~18.2]{shubin} for the corresponding result and expansion formula in the case of real-valued phase functions). The expansion coefficients are specified by (see \cite[equation (2.28)]{ms})
\begin{multline*} C(x,y,t) \sim \sum_\alpha \tilde e^{(\alpha)}(x, t\Phi'_x(x,y)) \frac{D_z^\alpha (A(z,y,t) e^{i \rho(z,x,y,t)})}{\alpha !} \Big{|}_{z = x}\\ \text{in }S^{n+m}_{1,0}(D \times D \times \RR_+),\end{multline*}
where $D^\alpha_z = (-i)^{|\alpha|} \partial^{\alpha}_z$ and $\tilde e^{(\alpha)} = \partial_\xi^{\alpha} \tilde e$ and $$\rho(z,x,y,t) = t\Phi(z,y) - t\Phi(x,y) - t(z-x) \cdot \Phi'_x(x,y).$$
Note that
\begin{equation*} \tilde e \sim \sum_{j \ge 0} \tilde e_j \text{ in } S^m_{\text{cl}}(D_{\CC}\times \CC^{2n+1}),\end{equation*}
where $D_{\CC} \subset \CC^{2n+1}$ is an open set such that $D_\CC \cap \RR^{2n+1} = D$. Also, $$\tilde e^{(\alpha)} \in S^{m-|\alpha|}_{\text{cl}}(D_\CC\times \CC^{2n+1}),$$ while $D^\alpha_z(A(z,y,t) e^{i \rho(z,x,y,t)})\big{|}_{z=x}$ is a classical symbol of order $n+\lfloor\frac{|\alpha|}2\rfloor$. Consequently, if $|\alpha| \ge 3$, then
\begin{equation*} \tilde e^{(\alpha)}(x, t\Phi'_x(x,y)) \frac{D^\alpha_z(A(z,y,t) e^{i \rho(z,x,y,t)})}{\alpha!}\Big{|}_{z=x} \in S^{n+m-2}_{\cl}(D \times D \times \RR_+).\end{equation*}

We now address the cases $|\alpha| = 0,1,2$. Note that $\rho(z,x,y,t)$ and its first derivatives (with respect to $z$) vanish when $z=x$. First, if $|\alpha| = 0$, we obtain
\begin{equation*} \tilde e^{(\alpha)}(x, t\Phi'_x(x,y)) \frac{D^\alpha_z(A(z,y,t) e^{i \rho(z,x,y,t)})}{\alpha!} \Big{|}_{z=x}= \tilde e(x, t\Phi'_x(x,y)) A(x,y,t),\end{equation*}
where (since $A_0(x,y) = \frac 1 {2\pi^{n+1}} + \mco(|(x,y)|^2)$)
\begin{multline*} \tilde e(x, t\Phi'_x(x,y)) A(x,y,t) \sim \tilde e_0(x, \Phi'_x(x,y)) A_0(x,y) t^{n+m} \\+ \Big{(}\tilde e_0(x, \Phi'_x(x,y)) A_1(x,y) + \tilde e_1(x, \Phi'_x(x,y)) A_0(x,y)\Big{)} t^{n+m-1} \\ \mod S^{n+m-2}_{\cl}(D \times D \times \RR_+)\\ \equiv \tilde e_0(x, \Phi'_x(x,y))A_0(x,y)t^{n+m} \\+ \left(\tilde e_0(x, \Phi'_x(x,y)) A_1(x,y) + \frac{\tilde e_1(x, \Phi'_x(x,y))}{2\pi^{n+1}}  + \mco(|(x,y)|^2)\right) t^{n+m-1}\\ \mod S^{n+m-2}_{\cl}(D \times D \times \RR_+). \end{multline*}
Next, consider the case $|\alpha| = 1$. Then
\begin{equation*} \tilde e^{(\alpha)}(x, t\Phi'_x(x,y)) \in S^{m-1}_{\cl}(D \times D \times \RR_+),\end{equation*}
and (since $A_0(x,y) = \frac 1 {2\pi^{n+1}} + \mco(|(x,y)|^2)$, and $\rho$ vanishes to second order when $z = x$)
\begin{multline*} D^\alpha_z(A(z,y,t) e^{i \rho(z,x,y,t)})\big{|}_{z=x} = D^\alpha_z(A(z,y,t))\big{|}_{z=x}\\ = \mco(|(x,y)|)t^n \mod S^{n-1}_{\cl}(D \times D \times \RR_+).\end{multline*}
Thus,
\begin{multline*} \tilde e^{(\alpha)}(x, t\Phi'_x(x,y)) \frac{D^\alpha_z(A(z,y,t) e^{i \rho(z,x,y,t)})}{\alpha!} \Big{|}_{z=x}  \\= \mco(|(x,y)|) t^{n+m-1} \mod S^{n+m-2}_{\cl}(D \times D \times \RR_+).\end{multline*}
Now, if $|\alpha| = 2$, then $D^\alpha_z = -\partial_{z_l, z_j}$ for some $1 \le j,l \le 2n+1$. Hence
\begin{equation*} D^\alpha_z(A(z,y,t) e^{i\rho(z,x,y,t)})\big{|}_{z=x} = - \partial_{x_j, x_l} A(x,y,t) -i A(x,y,t) \partial_{z_j, z_l} \rho(x,x,y,t).\end{equation*}
Also,
\begin{equation*} \partial_{z_j, z_l} \rho(x,x,y,t) = t \partial_{x_j, x_l} \Phi(x,y).\end{equation*}
Thus,
\begin{equation*} D^\alpha_z(A(z,y,t) e^{i\rho(z,x,y,t)})\big{|}_{z=x} = -i t^{n+1}A_0(x,y) \partial^\alpha_{x} \Phi(x,y) \mod S^{n}_{\cl}(D \times D \times \RR_+).\end{equation*}
Additionally,
\begin{equation*} \tilde e^{(\alpha)}(x, t \Phi'_x(x,y)) \sim t^{m-2} \partial_{\xi_j, \xi_l} \tilde e_0(x, \Phi'_x(x,y)) \mod S^{m-3}_{\cl}(D \times D \times \RR_+),\end{equation*}
so we obtain
\begin{multline*} \tilde e^{(\alpha)}(x, t\Phi'_x(x,y)) \frac{D^\alpha_z(A(z,y,t) e^{i \rho(z,x,y,t)})}{\alpha!}\Big{|}_{z=x} \\= - \frac i {\alpha!}  \partial^\alpha_{\xi} \tilde e_0(x, \Phi'_x(x,y)) A_0(x,y) \partial^\alpha_{x} \Phi(x,y)t^{n+m-1}\mod S^{n+m-2}_{\cl}(D \times D \times \RR_+)\\ = \left(-\frac i {2\alpha!\pi^{n+1}} \partial^\alpha_{\xi} \tilde e_0(x, \Phi'_x(x,y)) \partial^\alpha_{x} \Phi(x,y) + \mco(|(x,y)|^2) \right)t^{n+m-1} \\\mod S^{n+m-2}_{\cl}(D \times D \times \RR_+).\end{multline*}
Finally, since $A(x,y,t)$ is independent of $y_{2n+1}$, it is clear (in light of \cite[(2.28)]{ms}) that the same holds for $C_l(x,y)$ for any $l \ge 0$, hence without loss of generality also for $C(x,y,t)$.
\end{proof}
We will only require the values of $C_{1}$ at $x = y = 0$.
\begin{corollary}\label{c_1_at_0} In particular,
\begin{equation*} C_{1}(0,0) = \frac{R(0)e_0(0, -\omega_0(0)) +\sum_{m=1}^{2n} \partial^2_{\xi_m} e_0(0, -\omega_0(0))}{4 \pi^{n+1}}+ \frac{e_1(0, -\omega_0(0))}{2\pi^{n+1}} .\end{equation*}
Here, we view $e_0$, $e_1$ as functions on $T^*D$. We also note that since $\omega_0$ is real, it holds that $\tilde e_0(0, -\omega_0(0)) = e_0(0, -\omega_0(0))$ and $\tilde e_1(0, -\omega_0(0)) = e_1(0, -\omega_0(0))$. \end{corollary}
Now, we turn to study the relation (\ref{kern_comp}) between $T_E(x,y)$ and $\Pi(x,y)$, $Q_E(x,y)$. Recall that $dv_X(u) = \lambda(u) du$, as specified in Proposition \ref{canonic_coor}. Thus,
\begin{equation}\label{T_E_prod} \begin{aligned} &T_E(x,y) \equiv \int_D \Pi(x,u) Q_E(u,y) dv_X(u)\\ &= \int_D \int_0^\infty \int_0^\infty  e^{i s \Phi(x,u) + i t \Phi(u,y)} A(x,u,s) C(u,y,t)\lambda(u) ds dt du,\end{aligned}\end{equation}
hence we may apply Corollary \ref{toeplitz_exp} to $T_E(x,y)$.
\begin{proposition} There exists $B \in S^{n+m}_{\cl}(D \times D \times \RR_+)$, which is independent of $y_{2n+1}$, such that $$T_E(x,y) \equiv \int_0^\infty e^{it \Phi(x, y)} B(x,y,t) dt.$$
where
$$B(x,y,t) \sim \sum_{j \ge 0} B_{j}(x,y) t^{n+m-j}\ \text{in } S^{n+m}_{1,0}(D \times D \times \RR_+),$$
with
\begin{multline*} B_{0}(0,0) = \frac{e_0(0, -\omega_0(0))}{2\pi^{n+1}},\\ B_{1}(0,0) = \frac 1 {4\pi^{n+1}} \left[R(0) e_0(0, -\omega_0(0))-\square_b e_0(0, -\omega_0(0))\right] +\\ \frac 1 {4\pi^{n+1}} \left[ 2 e_1(0, -\omega_0(0))+ i\sum_{j=1}^{2n} \partial_{x_j, \xi_j} e_0(0, -\omega_0(0)) + \frac 1 2 \sum_{j=1}^{2n} \partial^2_{\xi_j} e_0(0, -\omega_0(0))\right].  \end{multline*} \end{proposition}
\begin{proof}
Write
\begin{equation*} \tilde e_{0,0}(x) = \tilde e_0(x, \Phi'_x(x,0)).\end{equation*} Applying Corollary \ref{toeplitz_exp} to $T_E(x,y)$ yields (in light of Propositions \ref{Q_E_prop}, \ref{hsiao_shen_indep} and Corollary \ref{c_1_at_0})
\begin{equation*} B_{0}(0,0) = \frac{e_0(0, -\omega_0(0))}{2\pi^{n+1}}\end{equation*}
and
\begin{multline*}
B_{1}(0,0) = \frac 1 {4\pi^{n+1}} \left[R(0)e_0(0, -\omega_0(0))-\square_b \tilde e_{0,0}(0)+2e_1(0, -\omega_0(0))\right]\\ +\frac 1 {4\pi^{n+1}}\sum_{m=1}^{2n} \partial^2_{\xi_m} e_0(0, -\omega_0(0))\end{multline*}
A straightforward computation produces
\begin{multline*} \square_b \tilde e_{0,0}(0) = \\-\frac 1 2 \sum_{j=1}^{2n}\left[ \partial^2_{x_j} e_0(0, -\omega_0(0)) +2i \partial^2_{x_j, \xi_j}e_0(0, -\omega_0(0)) -\partial^2_{\xi_j} e_0(0, -\omega_0(0))\right] \\-in \partial_{x_{2n+1}} e_0(0, -\omega_0(0))\\ = \square_b e_0(0, -\omega_0(0)) -i\sum_{j=1}^{2n} \partial_{x_j, \xi_j} e_0(0, -\omega_0(0)) +\frac 1 2 \sum_{j=1}^{2n} \partial_{\xi_j}^2 e_0(0, -\omega_0(0)).\end{multline*}
\end{proof}
We aim to express $B_{1}(0,0)$ in terms of globally defined objects. This will be achieved in steps, as follows. Let $$P : C^\infty(T^*X) \to C^\infty(T^*X)$$
be the differential operator specified in Corollary \ref{div_jhamvf}.
\begin{corollary}\label{intermediate_formula} Let $\mathcal E_0(x) =e_0(x, -\omega_0(x))$. Then
\begin{multline*} 4 \pi^{n+1}B_{1}(0,0) = R(0) \mathcal E_0(0) -\square_b \mathcal E_0(0) + P(e_0)(0, -\omega_0(0))\\+ 2 e_1(0, -\omega_0(0)) + i \sum_{j=1}^{2n} \partial_{x_j, \xi_j} e_0(0, -\omega_0(0)). \end{multline*} \end{corollary}
\begin{proof}
First, identifying $T^*D$ with $D \times \RR^{2n+1}$, we note that
\begin{equation*} \omega_0(x) = (x_2, -x_1, x_4, -x_3, ..., x_{2n}, -x_{2n-1}, 1+q_2(x)) + \mco(|x|^3),\end{equation*}
hence in light of the properties of $q_2(x)$ (specified in Proposition \ref{canonic_coor}), it suffices to consider
\begin{equation*} \hat \mce_0(x) = e_0(x, -\hat \omega_0(x)),\ \hat \omega_0(x) = (x_2, -x_1,..., x_{2n}, -x_{2n-1}, 1),\end{equation*}
since $\hat \mce_0(0) = \mce_0(0)$, $\square_b \hat \mce_0(0) = \square_b \mce_0(0)$. Next, for $j = 1,...,n$,
\begin{equation*} \partial_{x_{2j}} \hat \mce_0(x) = \partial_{x_{2j}} e_0(x, -\hat \omega_0(x)) - \partial_{\xi_{2j-1}} e_0(x, -\hat \omega_0(x)),\end{equation*}
and
\begin{equation*} \partial^2_{x_{2j}} \hat \mce_0(x) = \partial_{x_{2j}}^2 e_0(x, - \hat \omega_0(x)) - 2\partial_{\xi_{2j-1},x_{2j}} e_0(x, -\hat \omega_0(x))  + \partial^2_{\xi_{2j-1}} e_0(x, -\hat \omega_0(x)).\end{equation*}
Similarly,
\begin{equation*} \partial_{x_{2j-1}} \hat \mce_0(x) = \partial_{x_{2j-1}} e_0(x, -\hat \omega_0(x)) + \partial_{\xi_{2j}} e_0(x, -\hat \omega_0(x))\end{equation*}
and
\begin{equation*} \partial^2_{x_{2j-1}} \hat \mce_0(x) = \partial^2_{x_{2j-1}} e_0(x, -\hat \omega_0(x)) + 2\partial_{x_{2j-1}, \xi_{2j}} e_0(x, -\hat \omega_0(x)) + \partial^2_{\xi_{2j}} e_0(x, -\hat \omega_0(x)).\end{equation*}
Summing over $j = 1,...,n$, we see that
\begin{multline*} \sum_{j=1}^{2n} \partial^2_{x_j} \hat \mce_0(x) = \sum_{j=1}^{2n}\left( \partial^2_{x_j} e_0(x, -\hat \omega_0(x)) + \partial^2_{\xi_j} e_0(x, -\hat \omega_0(x))\right)\\ + 2\sum_{j=1}^n \left(\partial_{x_{2j-1}, \xi_{2j}} e_0(x, -\hat \omega_0(x)) - \partial_{x_{2j}, \xi_{2j-1}} e_0(x, -\hat \omega_0(x))\right).\end{multline*}
Noting that $\partial_{x_{2n+1}} \mce_0(0) = \partial_{x_{2n+1}} e_0(0, -\hat \omega_0(0))$, it follows that
\begin{multline*} -\square_b \mce_0(0) = -\square_b e_0(0, -\omega_0(0)) + \frac 1 2 \sum_{j=1}^{2n} \partial^2_{\xi_j} e_0(0, -\omega_0(0))\\ + \sum_{j=1}^n \left[\partial_{x_{2j-1}, \xi_{2j}} e_0(0, -\omega_0(0)) - \partial_{x_{2j}, \xi_{2j-1}} e_0(0, - \omega_0(0))\right].\end{multline*}
In light of Corollary \ref{div_jhamvf},$$\sum_{j=1}^n \left[\partial_{x_{2j-1}, \xi_{2j}} e_0(0, -\omega_0(0)) - \partial_{x_{2j}, \xi_{2j-1}} e_0(0, - \omega_0(0))\right] = -P(e_0)(0, -\omega_0(0)),$$
as required.
\end{proof}

We also require the following identity.
\begin{lemma}\label{princ_symb_id} The principal symbol $e_0 \in C^\infty(T^* X\setminus 0)$ of $E$ satisfies $$m \mct e_0(0, -\omega_0(0)) =\partial_{x_{2n+1}, \xi_{2n+1}} e_0(0, -\omega_0(0)).$$ \end{lemma}
\begin{proof}
Fix $(x,\xi) \in T^*D$ such that $|\xi| > 1$. Let $g(t) = e_0(x, t\xi)$, $t > 0$. Then by homogeneity, $$g'(t) = m t^{m-1} e_0(x, \xi).$$
On the other hand, using the chain rule, $$g'(t) = \sum_{j=1}^{2n+1} \partial_{\xi_j} e_0(x, t\xi) \xi_j.$$
Putting $t = 1$, we conclude that $$m e_0(x, \xi) = \sum_{j=1}^{2n+1} \partial_{\xi_j} e_0(x, \xi) \xi_j.$$
Differentiating this formula with respect to $x_{2n+1}$, we obtain $$m \partial_{x_{2n+1}} e_0(x, \xi) = \sum_{j=1}^{2n+1} \partial_{x_{2n+1}, \xi_j} e_0(x, \xi) \xi_j.$$
Thus for $x = 0$ and $\xi_\varepsilon = (0, -1-\varepsilon)$, $$m \partial_{x_{2n+1}} e_0(0, \xi_{\varepsilon}) = -(1+\varepsilon)\partial_{x_{2n+1}, \xi_{2n+1}} e_0(0, \xi_\varepsilon).$$
By continuity, we may set $\varepsilon = 0$ to conclude that $$m \partial_{x_{2n+1}} e_0(0,-\omega_0(0)) = -\partial_{x_{2n+1}, \xi_{2n+1}} e_0(0, -\omega_0(0)).$$
Since $\mct(0) = -\partial_{x_{2n+1}}$, this is as required.
\end{proof}
We now readily obtain the desired formula for $B_{1}(0,0)$.
\begin{corollary} Substituting the identity of Lemma \ref{princ_symb_id} into the formula for $B_{1}(0,0)$ of Corollary \ref{intermediate_formula}, and noting that $\mct e_0(0, -\omega_0(0)) = \mct \mce_0(0)$, we obtain
\begin{multline*} 4 \pi^{n+1}B_{1}(0,0) = R(0) \mathcal E_0(0) -\square_b \mathcal E_0(0)  + P(e_0)(0, -\omega_0(0))-i m \mct \mce_0(0)\\+ 2 e_1(0, -\omega_0(0)) + i \sum_{j=1}^{2n+1} \partial_{x_j, \xi_j} e_0(0, -\omega_0(0)). \end{multline*}
The subprincipal symbol $e_{\sub}$ of $E$ with respect to the $1$-density $|dv_X|$, as specified in Lemma \ref{density_subprinc_symb}, satisfies (in light of Proposition \ref{canonic_coor}) $$e_{\sub}(0, -\omega_0(0)) = e_1(0,-\omega_0(0)) + \frac i 2 \sum_{j=1}^{2n+1} \partial_{x_j, \xi_j} e_0(0, -\omega_0(0)).$$
Hence, we finally conclude that \begin{multline*}4\pi^{n+1} B_1(0,0) = \\ R(0) \mathcal E_0(0) - \square_b \mathcal E_0(0) + P(e_0)(0, -\omega_0(0)) - i m \mct \mce_0(0)+2 e_{\sub}(0, -\omega_0(0)).\end{multline*} \end{corollary}
\subsection{Expansions of products}
Let $(D, x)$ be a coordinate patch on $X$ as in Proposition \ref{canonic_coor}. Throughout, given a differential operator $Q : C^\infty(D) \to C^\infty(D)$, we denote by $$Q_x, Q_y : C^\infty(D \times D) \to C^\infty(D \times D)$$ the operators obtained by applying $Q$ on the $x$ and $y$ variables (respectively). Let $$\mca \in S^l_{\cl}(D \times D \times \RR_+),\ \mcb \in S^m_{\cl}(D \times D \times \RR_+),$$
where $l, m \in \RR$, and write
\begin{equation*} \begin{aligned} & \mca \sim \sum_{j \ge 0} \mca_j t^{l-j}\ \text{in } S^l_{1,0}(D \times D \times \RR_+),\\ & \mcb \sim \sum_{j \ge 0} \mcb_j t^{m-j}\ \text{in } S^m_{1,0}(D \times D \times \RR_+).\end{aligned} \end{equation*}
Assume that $\mca$, $\mcb$ and $\mca_j$, $\mcb_j$ (for all $j \ge 0$) are independent of $y_{2n+1}$.

Let $\Phi \in C^\infty(D \times D)$ denote the phase function of Proposition \ref{hsiao_shen_indep}. Define the oscillatory integral $$T(x,y) = \int_D \int_0^\infty \int_0^\infty e^{i(s \Phi(x,u) + t \Phi(u, y))} \mca(x,u, s) \mcb(u, y, t) \lambda(u) du ds dt.$$
The main result of the present subsection is as follows.
\begin{proposition}\label{product_expansion} There exists $\mcc \in S^{l+m-n}_{\cl}(D \times D \times \RR_+)$ such that $$T(x,y) \equiv \int_0^\infty e^{it \Phi(x,y)} \mcc(x,y,t)dt.$$
Moreover, $\mcc$ is independent of $y_{2n+1}$,
$$\mcc(x,y,t) \sim \sum_{j \ge 0} \mcc_j(x,y) t^{l+m-n-j}\ \text{in } S^{l+m-n}_{1,0}(D \times D \times \RR_+),$$
where $\mcc_j$, $j \ge 0$ are independent of $y_{2n+1}$, and
\begin{equation}\label{gen_C0C1}\begin{aligned} \mcc_0(0,0) = &\ 2\pi^{n+1} \mca_0(0,0) \mcb_0(0,0),\\ \frac{\mcc_1(0,0)}{\pi^{n+1}} = &- \mca_0(0,0) \mcb_0(0,0) R(0)\\ & + 2\big{(}\mca_0(0,0) \mcb_1(0,0) + \mca_1(0,0) \mcb_0(0,0) \big{)}\\ & -\big{(}\mca_0(0,0) \square_{b,x} \mcb_0(0,0) + \mcb_0(0,0) \square_{b,y} \mca_0(0,0)\big{)}\\ &+2i(n-l) \mca_0(0,0)\mct_x \mcb_0(0,0) + \sum_{j=1}^{2n} \partial_{y_j} \mca_0(0,0) \partial_{x_j} \mcb_0(0,0). \end{aligned}\end{equation}\end{proposition}
The following special case of Proposition \ref{product_expansion} is needed for the proof of Theorem \ref{main_thm_exp_sec}.
\begin{corollary}\label{toeplitz_exp} Assume that $\mca = A$, where $A \in S^n_{\cl}(D \times D \times \RR_+)$ is as specified in Proposition \ref{hsiao_shen_indep}. Then$$\mcc_0(0,0) = \mcb_0(0,0),$$
and $$\mcc_1(0,0) = \mcb_1(0,0)  -\frac 1 2 \square_{b,x} \mcb_0(0,0).$$ \end{corollary}
We establish Proposition \ref{product_expansion} in several steps.
\begin{lemma}\label{local_unique_noncanonic} Let $\Phi \in C^\infty(D \times D)$ be as in Proposition \ref{hsiao_shen_indep}. Then there exists $\mcc \in S^{l+m-n}_{\cl}(D \times D \times \RR_+)$, $$\mcc \sim \sum_{j \ge 0} \mcc_j t^{l+m-n-j}\ \text{in } S^{l+m-n}_{1,0}(D \times D \times \RR_+),$$ such that
\begin{equation}\label{product_kernel_local} T(x,y) \equiv \int_0^\infty e^{it \Phi(x,y)} \mcc(x,y,t) dt.\end{equation}
Additionally, we may assume without loss of generality that $\mcc$, $\mcc_j$ (for all $j \ge 0$) are independent of $y_{2n+1}$.\end{lemma}
\begin{proof}
Note that by the change of variable $s = \sigma t$,
\begin{equation}\label{for_diag_comp}T(x,y) = \int_D \int_0^\infty \int_0^\infty e^{it(\sigma \Phi(x,u) + \Phi(u, y))} t \mca(x, u, \sigma t) \mcb(u, y, t) \lambda(u) du d\sigma dt.\end{equation}
Let $\Psi(x,y,u,\sigma) = \sigma \Phi(x,u) + \Phi(u,y)$. Clearly, the fact that $\im \Phi \ge 0$ implies that $\im \Psi\ge 0$ (for $\sigma \ge 0$). Denote \begin{equation}\label{phase_0}\Psi_0(u, \sigma) = \Psi(0,0,u, \sigma).\end{equation} Then
\begin{equation*} \Psi_0(u, \sigma) = u_{2n+1}(\sigma-1) + (\sigma+1) \frac i 2 \sum_{j=1}^{2n} u_{j}^2 + \mco(|u|^4),\end{equation*}
hence $\Psi'_0(0,1) = 0$.
Additionally,
\begin{equation*} \Psi''_0(0,1) = \left(\begin{array}{ccc} 2i \Id_{2n} & 0 & 0 \\ 0 & 0 & 1 \\ 0 & 1 & 0 \end{array}\right).\end{equation*}
Thus, by the stationary phase formula (\cite[Theorem~2.3]{ms}),
\begin{equation}\label{T_E_ker} T(x,y) \equiv \int_0^\infty e^{i t \varphi(x,y)} \mcc(x,y,t) dt,\end{equation}
for some $\varphi \in C^\infty(D \times D)$ and $\mcc(x,y,t) \in S^{l+m-n}_{\cl}(D \times D \times \RR_+)$.

Note that $\varphi$ is determined regardless of $\mca$, $\mcb$. Specifically,
\begin{equation*} \varphi(x,y) = \tilde \Psi(x,y, \tilde U(x,y), \tilde \Sigma(x,y)),\end{equation*}
where $\tilde \Psi$ is an almost analytic extension of $$\Psi(x,y,u,\sigma) = \sigma \Phi(x,u) + \Phi(u,y),$$ and $\tilde U$, $\tilde \Sigma$ are determined by the equations (\cite[Lemma 2.1]{ms})
\begin{equation*} \begin{aligned} &\tilde \Phi(\tilde x,\tilde u) = 0,\\ &\tilde\sigma \partial_{\tilde y} \tilde \Phi(\tilde x,\tilde u) + \partial_{\tilde x} \tilde \Phi(\tilde u,\tilde y) = 0,\end{aligned} \end{equation*}
where now $\tilde \Phi$ is an almost analytic extension of $\Phi$, and $\tilde x$, $\tilde u$, $\tilde y$, $\tilde \sigma$ are complex. These equations do not depend on $\tilde y_{2n+1}$, hence the same holds for their solutions $\tilde U$, $\tilde \Sigma$, and it follows that
$$\tilde \Psi(x,y, \tilde U(x,y),\tilde \Sigma(x,y)) = y_{2n+1} + G(x,y')$$
for some smooth function $G$ (here $y' = (y_1, ..., y_{2n})$). Denote $$\Phi(x,y) = y_{2n+1} + g(x,y').$$Then as shown in \cite[Sect. 8]{hsiaomari2}, $$G(x,y') = g(x,y') + \mco(|x-y|^\infty).$$ Next (see \cite[(2.8)]{ms}), the fact that $\tilde U$, $\tilde \Sigma$ are independent of $\tilde y_{2n+1}$ implies that $\mcc(x,y,t)$ and $\mcc_j(x,y)$, $j \ge 0$, may be chosen (without loss of generality) to be independent of $y_{2n+1}$.
Finally, since $\varphi = \Phi + \mco(|x-y|^\infty)$, we conclude that $$T(x,y) \equiv \int_0^\infty e^{it (\Phi(x,y) + \mco(|x-y|^\infty))} \mcc(x,y,t) dt \equiv \int_0^\infty e^{it \Phi(x,y)} \mcc(x,y,t)dt.$$
\end{proof}
Next, we would like to establish formulas (\ref{gen_C0C1}) for $\mcc_0(0,0)$ and $\mcc_1(0,0)$.
\begin{lemma}\label{computation_lemma} Let $\mcc \in S^{l+m-n}_{\cl}(D \times D \times \RR_+)$ be as in Lemma \ref{local_unique_noncanonic}. Then \begin{equation*}\begin{aligned} \mcc(x,y,t) &\sim \sum_{j \ge 0} \mcc_j(x,y) t^{l+m-n -j}\ \text{in } S^{l+m-n}_{1,0}(D \times D \times \RR_+),\\ \mcc_0(0,0) = &\ 2\pi^{n+1} \mca_0(0,0) \mcb_0(0,0),\\ \frac{\mcc_1(0,0)}{\pi^{n+1}} = &- \mca_0(0,0) \mcb_0(0,0) R(0)\\ & + 2\big{(}\mca_0(0,0) \mcb_1(0,0) + \mca_1(0,0) \mcb_0(0,0) \big{)}\\ & -\big{(}\mca_0(0,0) \square_{b,x} \mcb_0(0,0) + \mcb_0(0,0) \square_{b,y} \mca_0(0,0)\big{)}\\ &+2i(n-l) \mca_0(0,0)\mct_x \mcb_0(0,0) + \sum_{j=1}^{2n} \partial_{y_j} \mca_0(0,0) \partial_{x_j} \mcb_0(0,0). \end{aligned}\end{equation*} \end{lemma}
\begin{proof}
Consider (\ref{for_diag_comp}) with $x = y= 0$. Let
\begin{equation*} I(t) = t\int_D \int_0^\infty e^{i t \Psi_0(u, \sigma)} \Gamma(u, \sigma, t) du d\sigma,\end{equation*}
where $\Psi_0$ is as in (\ref{phase_0}) and
\begin{equation*} \Gamma(u, \sigma, t) = \mca(0, u, \sigma t) \mcb(u, 0, t)\lambda(u).\end{equation*}
Then
\begin{equation*} I(t) \sim \mcc_0(0,0) t^{l+m-n} + \mcc_1(0,0) t^{l+m-n-1} + \mco(t^{l+m-n-2}).\end{equation*}
Introduce functions $\gamma_0$, $\gamma_1$ such that $$\Gamma(u, \sigma, t) = \gamma_0(u, \sigma) t^{l+m} + \gamma_1(u, \sigma)t^{l+m-1} + \mco(t^{l+m-2}).$$
Then (\cite[Theorem 7.7.5]{hormander1})
\begin{equation*} I(t) = \frac{ \gamma_0(0,1)t^{l+m+1} + \left(\gamma_1(0,1) + (L_1 \gamma_0)(0,1)\right)t^{l+m} + \mco(t^{l+m-1})}{ \left(\det\left(\frac{t \Psi''_0(0,1)}{2\pi i} \right) \right)^{\frac 1 2}}.\end{equation*}
Here, $L_1$ is the differential operator specified by
\begin{equation*} L_1 v = \frac 1 i \sum_{\mu=0}^2 \frac{\langle \Psi_0''(0,1)^{-1} \vec D, \vec D \rangle^{\mu+1} (h^\mu v)(0,1)}{\mu! (\mu+1)! 2^{\mu + 1}},\end{equation*}
where $\langle \cdot, \cdot \rangle$ is the complex bilinear dot product,
$$\vec D = -i(\partial_{u_1}, ..., \partial_{u_{2n+1}}, \partial_\sigma),$$ so that
\begin{equation*} \langle \Psi''_0(0,1)^{-1} \vec D, \vec D\rangle = -2 \partial_{u_{2n+1}} \partial_\sigma  + \frac i {2} \sum_{j=1}^{2n} \partial_{u_{j}}^2,\end{equation*}
and 
\begin{equation*} \begin{aligned} h(u, \sigma) &= \Psi_0(u, \sigma) - \frac 1 2 \langle \Psi''_0(0,1)(u, \sigma-1), (u, \sigma-1) \rangle\\ &= \frac i 2 (\sigma-1) \sum_{j=1}^{2n} u_{j}^2 + \mco(|u|^4).\end{aligned}\end{equation*}

First, note that
\begin{equation*} \det \left(\frac {t \Psi''(0,1)}{2\pi i} \right) = \frac{t^{2n+2}}{2^2 \pi^{2n+2}},\end{equation*}
hence
\begin{multline*} I(t) = 2 \pi^{n+1}\gamma_0(0,1)t^{l+m-n} + 2\pi^{n+1} \left(\gamma_1(0,1) + (L_1 \gamma_0)(0,1)\right)t^{l+m-n-1}\\ + \mco(t^{l+m-n-2}).\end{multline*}
Next,$$\begin{aligned} & \gamma_0(u, \sigma) = \mca_0(0, u) \mcb_0(u, 0) \lambda(u)\sigma^l,\\ & \gamma_1(0,1) = \mca_1(0,0) \mcb_0(0,0) + \mca_0(0,0) \mcb_1(0,0).\end{aligned}$$
We wish to compute $(L_1 \gamma_0)(0,1)$. Denote $$\gamma_0(u,\sigma) = \eta(u) \kappa(u, \sigma),$$
where $$\eta(u) = \mca_0(0, u) \mcb_0(u, 0),\ \kappa(u, \sigma) = \lambda(u) \sigma^l.$$
For $\mu = 0,1,2$, denote $$H_{\mu}(u,\sigma) = \left(h(u,\sigma)\right)^\mu \gamma_0(u, \sigma).$$
Also denote
\begin{equation*} h_0 =  (\sigma-1) \frac i 2 \sum_{j=1}^{2n} u_{j}^2,\ h_1 = h-h_0,\end{equation*}
and also
\begin{equation*}  \Delta = \sum_{j=1}^{2n} \partial_{u_{j}}^2.\end{equation*}
Consider the case $\mu = 2$. Then
\begin{multline*} \langle \Psi_0''(0,1)^{-1}\vec D, \vec D \rangle^3 =\\ -\frac i {8} \Delta^3 +\frac 3 2 \Delta^2 \partial_{\sigma,u_{2n+1}} +  6 i \Delta (\partial_{\sigma, u_{2n+1}})^2 -8 (\partial_{\sigma, u_{2n+1}})^3.\end{multline*}
Clearly, since $h(0,1) = \Delta h(0,1) = 0$, it holds that
\begin{equation*} \partial^3_{\sigma, u_{2n+1}} H_{2}(0,1) = \Delta \partial^2_{\sigma, u_{2n+1}}H_{2}(0,1) = 0.\end{equation*}
Also,
\begin{equation*} H_{2} = (h_0^2 + 2 h_0 h_1 + h_1^2) \gamma_0,\end{equation*}
where $$\Delta^2( h_0 h_1)(0,1)= \Delta^2(h_1^2)(0,1) = \partial_\sigma(h_0^2)(0,1)= 0,$$ hence
$\Delta^2 \partial_{\sigma, u_{2n+1}} H_{2}(0,1) = 0$. Similarly, writing $h_0 = (\sigma-1) \tilde h_0$,
$$\Delta^3 H_{2}= (\sigma-1) \Delta^3\left(((\sigma-1)\tilde h_0^2 + 2 \tilde h_0 h_1) \gamma_0\right) + \Delta^3\left(h_1^2  \gamma_0\right),$$
where $h_1^2 = \mco(|u|^8)$, and consequently $\Delta^3 H_{2}(0,1) = 0$. We conclude that
\begin{equation*} \langle \Psi''_0(0,1)^{-1} \vec D, \vec D \rangle^3 (H_{2})(0,1) = 0.\end{equation*}

Next, we consider the case $\mu = 1$, and note that
\begin{equation*} \langle \Psi''_0(0,1)^{-1} \vec D, \vec D \rangle^2 = -\frac 1 4 \Delta^2 -2 i \Delta \partial_{\sigma, u_{2n+1}} + 4 (\partial_{\sigma, u_{2n+1}})^2.\end{equation*}
First, $\partial_{\sigma, u_{2n+1}}^2 H_{1}(0,1) = 0$ since there is no differentiation with respect to $u_j$ for $j = 1,...,2n$. Next,
\begin{equation*} \Delta \partial_{\sigma, u_{2n+1}} H_{1} = \Delta \partial_{\sigma, u_{2n+1}}\left((\sigma-1) \tilde h_0 \gamma_0\right) + \Delta \partial_{\sigma, u_{2n+1}} \left(h_1 \gamma_0\right),\end{equation*}
and $\Delta \partial_{\sigma, u_{2n+1}} \left(h_1 \gamma_0 \right)(0,1) = 0$ since $h_1 = \mco(|u|^4)$. Also, since $\tilde h_0 = \mco(|u|^2)$,
\begin{equation*} \Delta \partial_{\sigma, u_{2n+1}}\left((\sigma-1) \tilde h_0 \gamma_0 \right)(0,1) = \partial_{u_{2n+1}}\gamma_0(0,1) \Delta \tilde h_0(u) = 2ni \partial_{u_{2n+1}}\gamma_0(0,1).\end{equation*}
Since $\kappa(0,1) = 1$ and $\partial_{u_{2n+1}} \kappa(0,1) = 0$ and $\mca_0$ is independent of $y_{2n+1}$,
\begin{multline*} \partial_{u_{2n+1}} \gamma_0(0,1) = \partial_{u_{2n+1}} \eta(0)\\ = \partial_{y_{2n+1}} \mca_0(0,0) \mcb_0(0,0) + \mca_0(0,0) \partial_{x_{2n+1}}\mcb_0(0,0) = \mca_0(0,0) \partial_{x_{2n+1}}\mcb_0(0,0).\end{multline*}
Thus,
\begin{equation*} \Delta \partial_{\sigma, u_{2n+1}} \left((\sigma-1) \tilde h_0 \gamma_0 \right)(0,1) = 2ni  \mca_0(0,0)\partial_{x_{2n+1}}\mcb_0(0,0).\end{equation*}
Finally, since $h_1 = \mco(|u|^4)$,
\begin{equation*} \Delta^2 H_{1}(0,1) = \gamma_0(0,1) \Delta^2 h_1(0,1) = \mca_0(0,0) \mcb_0(0,0) \Delta^2 h_1(0,1).\end{equation*}
We conclude that
\begin{multline*} \langle \Psi''_0(0,1)^{-1} \vec D, \vec D \rangle^2(H_{1})(0,1) =-\frac 1 4 \Delta^2 H_{1}(0,1) - 2 i \Delta \partial_{\sigma,u_{2n+1}} H_{1}\\ = -\frac 1 4 \mca_0(0,0) \mcb_0(0,0) \Delta^2 h_1(0,1) +4n \mca_0(0,0)\partial_{x_{2n+1}}\mcb_0(0,0).\end{multline*}
Finally, when $\mu = 0$, we recall that $\partial_{u_j} \lambda(0) = 0$ for $j = 1,...,2n+1$, so that
\begin{multline*} \partial_{\sigma, u_{2n+1}} (\gamma_0)(0,1)= l  \partial_{u_{2n+1}}\eta(0)\\ = l \left(\partial_{y_{2n+1}}\mca_0(0,0) \mcb_0(0,0) + \mca_0(0,0) \partial_{x_{2n+1}}\mcb_0(0,0) \right)\\ =l \mca_0(0,0) \partial_{x_{2n+1}}\mcb_0(0,0) ,\end{multline*}
and
\begin{equation*} \Delta(\gamma_0)(0,1) = \eta(0) \Delta\lambda(0) + \Delta \eta(0) = \mca_0(0,0) \mcb_0(0,0) \Delta\lambda(0) + \Delta \eta(0).\end{equation*}
Also,
\begin{equation*} \Delta \eta(0) = \mcb_0(0,0)(\Delta_y \mca_0)(0,0) + \mca_0(0,0)(\Delta_x \mcb_0)(0,0) + 2 \sum_{j=1}^{2n} \partial_{y_j} \mca_0(0,0) \partial_{x_j} \mcb_0(0,0).\end{equation*}
Thus,
\begin{multline*} \langle \Psi''_0(0,1)^{-1} \vec D, \vec D \rangle(H_{0})(0,1) =\frac i 2 \mca_0(0,0) \mcb_0(0,0) \Delta \lambda(0)\\ - 2l  \mca_0(0,0) \partial_{x_{2n+1}}\mcb_0(0,0) \\ + \frac i 2\left(  \mcb_0(0,0)(\Delta_y \mca_0)(0,0) + \mca_0(0,0)(\Delta_x \mcb_0)(0,0) + 2 \sum_{j=1}^{2n} \partial_{y_j} \mca_0(0,0) \partial_{x_j} \mcb_0(0,0)\right).\end{multline*}
Putting everything together, we obtain
\begin{multline*} iL_1(\gamma_0) =  \sum_{\mu = 0}^2 \frac{\langle \Psi''_0(0,1)^{-1} \vec D, \vec D \rangle^{\mu + 1}(H_{\mu})(0,1)}{\mu! (\mu+1)! 2^{\mu + 1}}\\ = \frac 1 8\left(-\frac 1 4 \mca_0(0,0) \mcb_0(0,0) \Delta^2 h_1(0,1) + 4n \mca_0(0,0)\partial_{x_{2n+1}} \mcb_0(0,0)\right)\\ + \frac 1 2 \left(\frac i 2 \mca_0(0,0) \mcb_0(0,0) \Delta \lambda(0) -2l\mca_0(0,0) \partial_{x_{2n+1}} \mcb_0(0,0)\right)\\ +\frac i 4 \left(\mcb_0(0,0) \Delta_y \mca_0(0,0) +\mca_0(0,0) \Delta_x \mcb_0(0,0) + 2 \sum_{j=1}^{2n} \partial_{y_j} \mca_0(0,0) \partial_{x_j} \mcb_0(0,0) \right)\\ = \mca_0(0,0) \mcb_0(0,0) \left(-\frac 1 {32} \Delta^2 h_1(0,1) + \frac i 4 \Delta \lambda(0) \right)\\ +\left(\frac{n-2l} 2 \right)  \mca_0(0,0) \partial_{x_{2n+1}} \mcb_0(0,0) \\+ \frac i 4 \left(\mcb_0(0,0) \Delta_y\mca_0(0,0) + \mca_0(0,0) \Delta_x \mcb_0(0,0)\right) \\+ \frac i 2 \sum_{j=1}^{2n} \partial_{y_j} \mca_0(0,0) \partial_{x_j} \mcb_0(0,0).\end{multline*}
Now, when $\mca = \mcb = A$, where $A\in S^n_{\cl}(D \times D \times \RR_+)$ is as in Proposition \ref{hsiao_shen_indep}, noting the properties of $A_0$, we see that $$i L_1 \gamma_0 =\frac 1 {4\pi^{2n+2}} \left(-\frac 1 {32} \Delta^2 h_1(0,1) + \frac i 4 \Delta \lambda(0) \right) = \frac {iL_1 \kappa(0,1)} {4\pi^{2n+2}} .$$
Define $\widehat \omega_0(x) = -d_x \Phi(x,x)$
and note that $\widehat \omega_0 = \nu \omega_0$, where $\nu = 1 - q_2 + \mco(|x|^3)$ is specified in Proposition~\ref{hsiao_shen_indep}. Define $\widehat \lambda\in C^\infty(D)$ by $$\frac 1 {n!} \left(\frac{-d\widehat \omega_0}2\right)^{\wedge n} \wedge \widehat \omega_0 = \widehat \lambda dx_1 \wedge ... \wedge dx_{2n+1},$$
and note that $$\widehat \lambda = \lambda -(n+1) q_2 + \mco(|x|^3).$$
Let $\widehat R$ denote the scalar curvature of the Tanaka-Webster connection associated with $\widehat \omega_0$. In \cite[Sect. 3.4]{hsiao_shen}, it is shown that $$-\frac 1 {32} \Delta^2 h_1(0,1) =\frac i 2 \widehat R(0),\ \frac i 4 \Delta \widehat \lambda(0) = - i \widehat R(0).$$
Noting that $\Delta \widehat \lambda(0) = \Delta \lambda(0)$, we obtain that $$i L_1 \kappa(0,1) = -\frac i 2 \widehat R(0).$$
Finally, it may be verified (for example, using \cite[eq. (5.56)]{drago_tomas}) that $\widehat R(0) = R(0)$.
This allows us to conclude that (no longer assuming $\mca = \mcb = A$)
\begin{multline*} L_1 \gamma_0(0,1) = -\frac 1 2 \mca_0(0,0) \mcb_0(0,0) R(0) \\ +\left(\frac{n-2l} {2i} \right)  \mca_0(0,0) \partial_{x_{2n+1}} \mcb_0(0,0) \\+ \frac 1 4 \left(\mcb_0(0,0) \Delta_y\mca_0(0,0) + \mca_0(0,0) \Delta_x \mcb_0(0,0)\right) \\+ \frac 1 2 \sum_{j=1}^{2n} \partial_{y_j} \mca_0(0,0) \partial_{x_j} \mcb_0(0,0). \end{multline*}
Consequently, \begin{multline*} \mcc_1(0,0) = 2\pi^{n+1}(\gamma_1(0,1) + L_1 \gamma_0(0,1))\\ = 2\pi^{n+1}\left(\mca_1(0,0) \mcb_0(0,0) + \mca_0(0,0) \mcb_1(0,0)\right)\\ -\pi^{n+1} \mca_0(0,0) \mcb_0(0,0) R(0) \\ +i(2l-n)\pi^{n+1}  \mca_0(0,0) \partial_{x_{2n+1}} \mcb_0(0,0) \\+ \frac {\pi^{n+1}} 2 \left(\mcb_0(0,0) \Delta_y\mca_0(0,0) + \mca_0(0,0) \Delta_x \mcb_0(0,0)\right) \\+ \pi^{n+1} \sum_{j=1}^{2n} \partial_{y_j} \mca_0(0,0) \partial_{x_j} \mcb_0(0,0). \end{multline*}
Noting that\begin{equation*} i(2l-n) \partial_{x_{2n+1}} \mcb_0(0,0) + \frac 1 2 \Delta_x \mcb_0(0,0) = 2i(n-l) \mct_x \mcb_0(0,0) - \square_{b,x} \mcb_0(0,0),\end{equation*}
and that \begin{equation*} \frac 1 2 \Delta_y \mca_0(0,0) = -\square_{b,y} \mca_0(0,0),\end{equation*} we finally conclude that
\begin{multline*} \frac 1{\pi^{n+1}}\mcc_1(0,0) = 2 \left(\mca_1(0,0) \mcb_0(0,0) + \mca_0(0,0) \mcb_1(0,0)\right) -\mca_0(0,0) \mcb_0(0,0) R(0) \\ +\mca_0(0,0)(2i(n-l) \mct_x \mcb_0(0,0) - \square_{b,x} \mcb_0(0,0)) -\mcb_0(0,0) \square_{b,y} \mca_0(0,0)\\
+ \sum_{j=1}^{2n} \partial_{y_j} \mca_0(0,0) \partial_{x_j} \mcb_0(0,0). \end{multline*}
\end{proof}
\section{Distribution kernels represented by oscillatory integrals}\label{aux_sect}
In the present section, we specify certain results concerning the representations of distribution kernels as oscillatory integrals. Let $\Gamma$ denote the gamma function, let $\ZZ_-$ denote the negative integers, and let $$\gamma = \lim_{m \to \infty} \left(\sum_{j=1}^m \frac 1 j - \log m \right)$$ denote Euler's constant. The results of this section rely on the following classical formulas\footnote{The integral on the left hand side is a so-called ``finite part integral'', cf. \cite[Ch. 3.2]{hormander1}.}, which are valid for $x \ne 0$, $\re(x) \ge 0$,
\begin{equation}\label{classical_eq}\int_0^\infty e^{-xt} t^m dt = \left\{\begin{array}{ll} \Gamma(m+1) x^{-m-1} & \text{if } m \in \RR \setminus \ZZ_-,\\ \frac{(-1)^m}{(-m-1)!} x^{-m-1} \left(\log x + \gamma -\sum_{j=1}^{-m-1} \frac 1 j \right)& \text{if } m \in \ZZ_-.\end{array}\right.\end{equation}
The following lemma is useful when dealing with the ambiguity in the choice of phase function. 
\begin{lemma}[cf. {\cite[Lemma 3.1]{hsiao_shen}}]\label{phase_factorization} Let $D \subset \RR^n$ be a sufficiently small open set with $0 \in D$. Assume that $$F \in C^\infty(D),\ F(0) = 0,\ \im F \ge 0,\ dF \ne 0\ \text{if } \im F = 0.$$ Assume that $$G \in C^\infty(D),\ \re(G)(0) > 0,\ \im(FG) \ge 0,\ d(FG) \ne 0\ \text{if } \im(FG) =0.$$ Let $m \in \RR$. Then in the sense of oscillatory integrals, $$\int_0^\infty e^{it G(x) F(x)} t^m dt \equiv \int_0^\infty e^{it F(x)} \frac{t^m}{(G(x))^{m+1}} dt \mod C^\infty(D).$$ \end{lemma}
\begin{proof}
The case $m \in \ZZ$, $m \ge 0$ was proven in \cite[Lemma 3.1]{hsiao_shen}, and when $m \in \RR \setminus \ZZ_-$ the proof is the same. Namely, in the sense of distributions, $$\begin{aligned}\int_0^\infty e^{it G(x) F(x)} t^m dt &= \lim_{\varepsilon \to 0^+} \int_0^\infty e^{it G(x) F(x) - \varepsilon t}t^m dt \\ &= \lim_{\varepsilon \to 0^+} \frac {\Gamma(m+1)}{(-i G(x) F(x) +\varepsilon)^{m+1}} \\ &= \frac 1 {(G(x))^{m+1}} \lim_{\varepsilon \to 0^+} \frac{\Gamma(m+1)}{\left(-iF(x) + \frac \varepsilon {G(x)}\right)^{m+1}} \\ &=\frac 1 {(G(x))^{m+1}} \lim_{\varepsilon \to 0^+} \int_0^\infty e^{it F(x) -\frac{\varepsilon t}{G(x)}} t^m dt \\ &= \int_0^\infty e^{it F(x)} \frac{t^m}{(G(x))^{m+1}} dt.\end{aligned}$$
If $m \in \ZZ_-$ then writing $c_m = \frac{(-1)^m}{(-m-1)!}$ and $\gamma_m = \gamma -\sum_{j=1}^{-m-1} \frac 1 j$,\begin{multline*} \int_0^\infty e^{it G(x) F(x)} t^m dt =\\ \lim_{\varepsilon \to 0^+}c_m(-iG(x) F(x) + \varepsilon)^{-m-1} \left(\log(-i F(x) G(x) +\varepsilon) + \gamma_m \right) \equiv \\  \frac 1 {(G(x))^{m+1}} \lim_{\varepsilon \to 0^+} c_m\left(-iF(x) + \frac \varepsilon{G(x)}\right)^{-m-1} \left(\log\left(-i F(x) + \frac \varepsilon{G(x)}\right) + \gamma_m \right)\\ = \int_0^\infty e^{it F(x)} \frac{t^m}{(G(x))^{m+1}} dt.\end{multline*}
\end{proof}
The formulas (\ref{classical_eq}) lead to an alternative representation of the distribution kernels of Toeplitz operators, which is of the same type as the standard representation of the Szeg\H{o} kernel (\cite[Corollaire~1.7]{bdms}). In what follows, we use the standard notation $$(-i(\phi(x,y) + i0))^{-(n+m+1)} := \lim_{\varepsilon \to 0^+}(-i(\phi(x,y) + i\varepsilon))^{-(n+m+1)},$$
where the limit is in the sense of distributions; similarly for \mbox{$\log(-i(\phi(x,y)+i0))$}.
\begin{corollary} Let $m \in \RR$ and $E \in L^m_{\cl}(X)$. Let $(D,x)$ be a coordinate patch such that $$T_E(x,y) \equiv \int_0^\infty e^{it \phi(x,y)}b(x,y,t) dt \mod C^\infty(D \times D),$$
where $\phi$ is as specified in Theorem \ref{bdms_thm}, and $b \in S^{n+m}_{\cl}(D \times D \times \RR_+)$, $$b \sim \sum_{j \ge 0} b_{j} t^{n+m-j}\ \text{in } S^{n+m}_{1,0}(D \times D \times \RR_+).$$ If $m \not \in \ZZ$, then there exists $F \in C^\infty(D \times D)$ which satisfies\footnote{Here, the symbol $\sim$ indicates that the left hand side and the right hand side have the same Taylor expansion at the diagonal.} $$F(x,y) \sim \sum_{j\ge 0}\Gamma(n+m+1-j) b_j(x,y) (-i \phi(x,y))^j,$$
such that $$T_E(x,y) \equiv \frac {F(x,y)} {(-i(\phi(x,y) + i0))^{n+m+1}} \mod C^\infty(D \times D).$$
If $m \in \ZZ$, $n+m < 0$, then there exists $G \in C^\infty(D \times D)$, $$G \sim \sum_{j \ge 0} \frac{(-1)^{n+m-j}}{(j-m-n-1)!} b_j(x,y)(-i\phi(x,y))^{j-m-n-1},$$
such that $$T_E(x,y) \equiv G(x,y) \log(-i(\phi(x,y) + i0)) \mod C^\infty(D \times D).$$
If $m \in \ZZ$, $n+m \ge 0$, then there exist $F, G \in C^\infty(D \times D)$, $$F = \sum_{j=0}^{n+m} (n+m-j)!b_j(x,y) (-i \phi(x,y))^j \mod \phi^{n+m+1},$$
and $$G \sim \sum_{j \ge 0} \frac{(-1)^{j+1}}{j!} b_{n+m+1+j}(x,y)(-i\phi(x,y))^j,$$
such that $$T_E(x,y) \equiv \frac{F(x,y)}{(-i(\phi(x,y)+i0))^{n+m+1}} + G(x,y) \log(-i(\phi(x,y)+i0)).$$
\end{corollary}
The following uniqueness result is central to the proof of Theorem \ref{main_thm}.
\begin{lemma}[cf. {\cite[Lemma 3.2]{hsiao_shen}}]\label{pointwise_uniqueness} Let $p \in X$ and let $(D, x)$ be a coordinate patch with $p \in D$. Assume that $\phi_1, \phi_2 \in C^\infty(D \times D)$ satisfy (\ref{phi_props}), and that they are both equivalent, in the sense of \cite{ms}, to the phase function $\phi$ of Theorem \ref{bdms_thm}. Assume further that $$\mct^2_y \phi_1(p,p) = \mct^2_y \phi_2(p,p) = 0.$$
Let $\alpha, \beta \in S^{n+m}_{\cl}(D \times D \times \RR_+)$, where $m \in \RR$, such that $$\begin{aligned} &\alpha \sim \sum_{j \ge 0} \alpha_j t^{n+m-j} \ \text{in } S^{n+m}_{1,0}(D \times D \times \RR_+),\\ &\beta \sim \sum_{j \ge 0} \beta_j t^{n+m-j} \ \text{in } S^{n+m}_{1,0}(D \times D \times \RR_+).\end{aligned}$$
Assume that $E \in L^m_{\cl}(X)$ satisfies $$\begin{aligned} & T_E(x,y) \equiv \int_0^\infty e^{it \phi_1(x,y)} \alpha(x,y,t) dt \equiv \int_0^\infty e^{it \phi_2(x,y)} \beta(x,y,t) dt\\ & \alpha_0(p,p) = \beta_0(p,p),\\ &\mct_y \alpha_0(p,p) = \mct_y \beta_0(p,p) = 0. \end{aligned}$$
Then $\alpha_1(p,p) = \beta_1(p,p)$.
\end{lemma}
\begin{proof}
We consider coordinates $x = (x_1, ...,x_{2n+1})$ on $D$ such that $x(p) = 0$ and $\mct = -\partial_{x_{2n+1}}$.
Then, we note that (\cite[proof of Lemma 3.2]{hsiao_shen}) there exists $f \in C^\infty(D \times D)$ such that $$\begin{aligned} &\phi_2(x,y) = f(x,y) \phi_1(x,y) + \mco(|x-y|^\infty),\\ &f(x,x) = 1,\ \mct_y f(0,0) = 0.\end{aligned}$$
so we may assume without loss of generality that $\phi_2 = f \phi_1$. 

If $m \in \RR \setminus \ZZ$, then \begin{multline*}\frac{\Gamma(n+m+1) \alpha_0 + \Gamma(n+m)\alpha_1(-i\phi_1) + \mco(|x-y|^2)}{(-i(\phi_1 + i0))^{n+m+1}} \\ \equiv \frac{\Gamma(n+m+1) \beta_0 + \Gamma(n+m)\beta_1(-if\phi_1) + \mco(|x-y|^2)}{(-i(f\phi_1 + i0))^{n+m+1}}.\end{multline*}
Thus, when $x = 0$, $y = (0, y_{2n+1})$, \begin{multline*}f^{n+m+1}(\Gamma(n+m+1) \alpha_0+\Gamma(n+m) \alpha_1(-i\phi_1) + \mco(|y_{2n+1}|^2)) \\= \Gamma(n+m+1) \beta_0 + \Gamma(n+m) \beta_1(-i f \phi_1)+\mco(|y_{2n+1}|^2)+(-i f \phi_1)^{n+m+1} S\end{multline*}
for some smooth function $S$. Since $f(0,(0,y_{2n+1})) = 1 + \mco(|y_{2n+1}|^2)$,\begin{multline*}\Gamma(n+m+1) \alpha_0 + \Gamma(n+m) \alpha_1(-i\phi_1) = \Gamma(n+m+1) \beta_0 +\Gamma(n+m) \beta_1 (-i \phi_1) \\+ \mco(|y_{2n+1}|^2)+(-i \phi_1)^{n+m+1} S.\end{multline*}
Since (noting that $\alpha_0(0,0) = \beta_0(0,0)$) $$\begin{aligned} &\alpha_0(0,(0,y_{2n+1})) = \alpha_0(0,0) + \mco(|y_{2n+1}|^2),\\ &\beta_0(0,(0,y_{2n+1})) = \alpha_0(0,0) + \mco(|y_{2n+1}|^2),\end{aligned}$$
we obtain $$\Gamma(n+m)(-i \phi_1)(\alpha_1- \beta_1) + \mco(|y_{2n+1}|^2)=(-i \phi_1)^{n+m+1}S.$$
Now, note that $\phi_1(0,(0,y_{2n+1})) = y_{2n+1} + \mco(|y_{2n+1}|^2)$. Thus, if $n+m+1 > 1$, we readily obtain that $\alpha_1(0,0) = \beta_1(0,0)$ as required. If $n+m+1 < 1$ (note that $m \not \in \ZZ$), then $$\Gamma(n+m)(-i \phi_1)^{-(n+m)}(\alpha_1 - \beta_1) + \mco(|y_{2n+1}|^{1-n-m}) = S,$$
which also implies (since $S$ is smooth) that $\alpha_1(0,0) = \beta_1(0,0)$.

Next, note that the case $n+m \in \ZZ$, $n +m > 0$ is essentially proven in \cite[Lemma 3.2]{hsiao_shen}. If $n+m = 0$, we note that $$\log(-i(f \phi_1 + i0)) \equiv \log(-i(\phi_1 + i 0)),$$ therefore \begin{multline*}\frac{\alpha_0}{(-i (\phi_1 + i0))} + (-\alpha_1 + \mco(|x-y|)) \log(-i(\phi_1+i0)) \equiv\\ \frac{\beta_0}{(-i(f \phi_1+i0))} +(-\beta_1 + \mco(|x-y|))\log(-i(\phi_1 + i0)).\end{multline*}
In particular, when $x = 0$ and $y = (0, y_{2n+1})$, using that $$\alpha_0(0,(0,y_{2n+1})) - \beta_0(0,(0,y_{2n+1})) = \mco(|y_{2n+1}|^2)$$ and $f(0,(0,y_{2n+1})) = 1+\mco(|y_{2n+1}|^2)$, we find that $$(\beta_1-\alpha_1 + \mco(|y_{2n+1}|))\log(-i \phi_1) = S(y_{2n+1})+\mco(|y_{2n+1}|)$$
for some smooth function $S$, which implies that $\alpha_1(0,0) = \beta_1(0,0)$.

Finally, if $n+m \in \ZZ$, $n+m < 0$, then \begin{multline*} (-i\phi_1)^{-n-m-1} \left(\frac{\alpha_0}{(-n-m-1)!}  - \frac{\alpha_1(-i\phi_1)}{(-n-m)!}  + \mco(|x-y|^2)\right)\log(-i(\phi_1+i0))  \equiv \\ (-if \phi_1)^{-n-m-1}\left(\frac{\beta_0}{(-n-m-1)!} - \frac{\beta_1(-if \phi_1)}{(-n-m)!} + \mco(|x-y|^2)\right)\log(-i(\phi_1+i0)).\end{multline*}
Again setting $x = 0$, $y = (0, y_{2n+1})$, we obtain $$(-i\phi_1)^{-n-m}\left(\beta_1 - \alpha_1\right)\log(-i \phi_1) =S(y_{2n+1})+ \mco(|y_{2n+1}|^{-n-m})$$
for some smooth function $S$, which readily implies that $\alpha_1(0,0) = \beta_1(0,0)$.
\end{proof}
\begin{corollary}[cf. {\cite[Lemma 1.1]{hsiao_shen}}]\label{uniqueness_corollary} Let $(D,x)$ be an open coordinate patch. Assume that $\phi_1, \phi_2 \in C^\infty(D \times D)$ satisfy (\ref{phi_props}), and that they are both equivalent, in the sense of \cite{ms}, to the phase function $\phi$ of Theorem \ref{bdms_thm}. Assume further that $$\mct^2_y \phi_1(x,x) = \mct^2_y \phi_2(x,x) = 0\text{ for all }x\in D.$$
Let $\alpha, \beta \in S^{n+m}_{\cl}(D \times D \times \RR_+)$, where $m \in \RR$, such that $$\begin{aligned} &\alpha \sim \sum_{j \ge 0} \alpha_j t^{n+m-j}\ \text{in } S^{n+m}_{1,0}(D \times D \times \RR_+),\\ &\beta \sim \sum_{j \ge 0} \beta_j t^{n+m-j}\ \text{in } S^{n+m}_{1,0}(D \times D \times \RR_+). \end{aligned}$$
Assume that $E\in L^m_{\cl}(X)$ satisfies $$\begin{aligned} &T_E(x,y) \equiv \int_0^\infty e^{it \phi_1(x,y)} \alpha(x,y,t) dt \equiv \int_0^\infty e^{it \phi_2(x,y)} \beta(x,y,t) dt,\\ &\alpha_0(x,x) = \beta_0(x,x)\ \text{for all }x \in D,\\ &\mct_y \alpha_0(x,x) = \mct_y \beta_0(x,x) = 0\ \text{for all }x \in D. \end{aligned}$$
Then $\alpha_1(x,x) = \beta_1(x,x)$ for all $x \in D$.\end{corollary}
\section{Distribution kernels of Toeplitz operators on CR orbifolds}\label{s-gue251129yyd}

In this section, we will establish asymptotic expansions for the Toeplitz operators on a compact not necessarily strictly pseudoconvex CR orbifold. We first recall the definition of CR orbifolds (see~\cite[Sect. 5]{khh}). 

\begin{definition}\label{d-gue181015}
Let $X$ be a Hausdorff topological space. 
We say that $X$ is a CR orbifold of dimension $2n+1$ with CR codimension $1$ if there exists a cover $U_i$ of $X$, which is closed under finite intersections, such that
\begin{itemize}
\item For each $U_i$, there exists a CR manifold $V_i$ of dimension $2n+1$ with CR codimension $1$, a finite group $\Gamma_i$ acting on $V_i$ with CR automorphisms and a $\Gamma_i$-invariant map $\Psi_i \colon V_i \rightarrow U_i$, which induces a homeomorphism $V_i / \Gamma_i \rightarrow U_i$. 
\item For each inclusion $U_i \subset U_j$, we have an injective group morphism $\varphi_{ij} \colon \Gamma_i \rightarrow \Gamma_j$ and a CR isomorphism $\Phi_{ij} \colon V_i \rightarrow \Psi_j^{-1}(U_i)$, which satisfies $\Phi_{ij}(g\cdot x) = \varphi_{ij} (g)\cdot\Phi_{ij} (x)$ for $x \in V_i$, $g \in \Gamma_i$ and fulfills $\Psi_j \circ \Phi_{ij} = \Psi_i$.
\end{itemize}
We call the tuple $(U_i, V_i, \Gamma_i, \Psi_i)$ an orbifold chart and the cover $U_i$ an orbifold atlas.

An orbifold is called effective if for all charts, the action of $\Gamma_i$ on $V_i$ is effective.

For $U \subset X$ open, a function $f \colon U \rightarrow \mathbb C$ is called a CR function if every lift of $f$ into a chart is a CR function. 

For an orbifold chart $(U, V, \Gamma, \Psi)$, we say that $X$ is strictly pseudoconvex on $U$ if $V$ is a strictly pseudoconvex CR manifold. 
\end{definition}

For any orbifold chart $(U, V, \Gamma, \Psi)$, we will always identify $U$ with $V/\Gamma$. 

From now on, we let $X$ be a compact orientable CR orbifold of dimension $2n+1$ with CR codimension one. 
We fix a Hermitian metric $\langle\,\cdot\,|\,\cdot\,\rangle$ on $TX\otimes_{\mathbb R}\mathbb C$ and we let $dv_X$ be the volume form on $X$ induced by $\langle\,\cdot\,|\,\cdot\,\rangle$ and let $(\,\cdot\,|\,\cdot\,)$ be the $L^2$ inner product on $C^\infty(X)$ induced by $dv_X$ and as in the smooth case, let $L^2(X)=L^2(X,dv_X)$ be the completion of $C^\infty(X)$ induced by $(\,\cdot\,|\,\cdot\,)$. All the standard notations, terminology and set up as in Section~\ref{prelim_sect} can be generalized to the orbifolds setting. For simplicity, we omit the details and refer the reader to~\cite{gh3}. We will use the same notations as in the smooth case. 

For $\lambda\geq0$, let 
\[\Pi_{\leq\lambda}:=1_{[0,\lambda]}(\square_b): L^2(X)\rightarrow L^2(X),\]
where $1_{[0,\lambda]}(\square_b)$ denotes the functional calculus of $\square_b$ with respect to $1_{[0,\lambda]}$. 
Let $\Pi_{\leq\lambda}(x,y)\in\mcd'(X \times X)$ be the distribution kernel of $\Pi_{\leq\lambda}$. For $\lambda=0$, we write $\Pi:=\Pi_{\leq\lambda}$, $\Pi(x,y)=\Pi_{\leq\lambda}(x,y)$.
We have  

\begin{theorem}\label{t-gue251130yyd}
With the notations used above, let $(U, V, \Gamma, \Psi)$
be an orbifold chart of $X$. Assume that $X$ is strictly pseudoconvex on $U$. Fix $\lambda>0$. We have 
		\begin{equation}\label{e-gue251130yyda}
        \Pi_{\leq\lambda}(x,y)=\frac{1}{|\Gamma|}\sum_{h,g\in\Gamma}\Bigl(\int _0^{+\infty}e^{it\, \phi(h\cdot\tilde{x},g\cdot \tilde y)}a(h\cdot\tilde{x},g\cdot \tilde{y},t)dt+F(h\cdot\tilde x,g\cdot\tilde y)\Bigr)\ \ \mbox{on $U\times U$} \,,\end{equation}
		where $\pi(\widetilde x)=x$, $\pi(\widetilde y)=y$, $\pi: V\To U$ is the natural projection, $F$ is a  smoothing operator on $V$ and 
		\[a( \widetilde{x}, \widetilde{y},t)\sim\sum^\infty_{j=0}a_j( \widetilde{x},  \widetilde{y})\,t^{n-j}\]
		in $S^{n}_{1, 0}(V\times V\times\mathbb{R}_+)$,
        $a_j( \widetilde{x},  \widetilde{y})\in C^\infty(V\times V)$. $j=0,1,\ldots$, 
        \begin{equation}\label{e-gue251210ycda}
        a_0(\widetilde{x},  \widetilde{x})=\frac{1}{2}\pi^{-n-1}|\det\mathcal{L}_{\widetilde{x}}|,\ \ \widetilde{x}\in V,\end{equation}
        where $\det\mathcal{L}_{\widetilde{x}}=\mu_1(\widetilde x)\cdots\mu_n(\widetilde x)$, $\mu_j(\widetilde x)$, $j=1,\ldots,n$, are the eigenvalues of 
        $\mathcal{L}_{\widetilde{x}}$ with respect to $\langle\,\cdot\,|\,\cdot\,\rangle$, 
		 the phase function $\phi$ is the same as in Theorem~\ref{bdms_thm},
  $a$ and the phase functions $\phi$ can be taken to be $\Gamma$-invariant which means 
  $a(g\cdot\tilde{x},g\cdot\tilde{y},t)=a(\tilde{x},\tilde{y},t)$, $\phi(g\cdot\tilde{x},\,g\cdot\tilde{y})=\phi(\tilde{x},\,\tilde{y})$, for every $g\in\Gamma$ and every $\tilde x, \tilde y \in V$.
\end{theorem}

\begin{proof}
From the construction in~\cite[part I, Chapter 8]{hsiao1}, there exist properly supported continuous operators 
		$A, S: C^\infty(V)\To C^\infty(V)$  such that
		\begin{equation}\label{e-gue140205II}
			\begin{split}
				&\mbox{$\square_bA+S=I$ on $V$},\\
                &\mbox{$A^*\square_b+S=I$ on $V$},\\
				&\square_bS\equiv0\ \ \mbox{on $V$},\\
				&S\equiv S^*\equiv S^2\ \ \mbox{on $V$},
			\end{split}
		\end{equation}
		where $S^*$ and $A^*$ are the formal adjoints of $S$ and $A$
		with respect to $(\,\cdot\,|\,\cdot\,)$ respectively and $S(\widetilde{x},\widetilde{y})$ satisfies
		\begin{equation}\label{e-gue240309yyd}
  S(\widetilde{x}, \widetilde{y})\equiv\int^{\infty}_{0}e^{i\phi(\widetilde{x}, \widetilde{y})t}s(\widetilde{x}, \widetilde{y}, t)\,\mathrm{d}t\ \ \mbox{on $V$}\end{equation}
		with a symbol $s\in S^{n}_{{\rm cl\,}}\big(V\times V\times\mathbb{R}_+\big)$  such that
		\begin{equation}  \label{e-gue140205III}\begin{split}
				&s(\widetilde{x}, \widetilde{y}, t)\sim\sum^\infty_{j=0}s_j(\widetilde{x}, \widetilde{y})t^{n-j}\quad\text{ in }S^{n}_{1, 0}
				\big(V\times V\times\mathbb{R}_+\big)\,,\\
				&s_j\in C^\infty\big(V\times V\big),\ \ j\in\mathbb N\cup\{0\},
		\end{split}\end{equation}
        $s_0$ satisfies \eqref{e-gue251210ycda} 
		and the phase function $\phi$ is the same as in Theorem~\ref{bdms_thm}.  

        Define on $U$ the following kernels
	\[
		\hat{S}({x},{y})=\frac{1}{|\Gamma|}\sum_{h,g\in\Gamma}S(h\cdot\widetilde{x},g\cdot \widetilde{y})\in\mcd'(U\times U),\]
	and
	\[\hat{A}(x,y)=\frac{1}{|\Gamma|}\sum_{g,h\in\Gamma} A( h\cdot\widetilde{x},g\cdot \widetilde{y})\in\mcd'(U\times U). \]
	 For every $u \in C^\infty_0(U)$, we define
	\[\hat{A}u=\frac{1}{|\Gamma|^2}\sum_{g,h\in\Gamma}\int_{V} A( h\cdot\widetilde{x},g\cdot \widetilde{y})u(\widetilde{y})dv_X\]
	which is in $C^\infty_0(U)$ and similarly one defines 
    \[\hat{S}u=\frac{1}{|\Gamma|^2}\sum_{g,h\in\Gamma}\int_{V} S( h\cdot\widetilde{x},g\cdot \widetilde{y})u(\widetilde{y})dv_X.\]
    Then, $\hat A$ and $\hat S$ are properly supported continuous operators:
    \[\hat A, \hat S: C^\infty(U)\To C^\infty(U).\]
    
   From \eqref{e-gue140205II}, we see that 
 \begin{equation}\label{e-gue251206yyd}
			\begin{split}
				&\mbox{$\square_b\hat{A}+\hat{S}=I$ on $U$},\\
                &\mbox{$\hat{A}^*\square_b+\hat{S}^*=I$ on $U$},\\
				&\square_b\hat{S}\equiv0\ \ \mbox{on $U$},\\
				&\hat{S}\equiv \hat{S}^*\equiv \hat{S}^2\ \ \mbox{on $U$},
			\end{split}
		\end{equation}
		where $\hat{S}^*$ and $\hat{A}^*$ are the formal adjoints of $\hat{S}$ and $\hat{A}$
		with respect to $(\,\cdot\,|\,\cdot\,)$ respectively. 

From \eqref{e-gue251206yyd}, we can repeat the argument in the proof of~\cite[Theorem 1.5]{hsiaomari2} and deduce that 
\[\Pi_{\leq\lambda}\equiv\hat{S}\ \ \mbox{on $U$}.\]
\end{proof}

We now introduce pseudodifferential operators on orbifolds. 

\begin{definition}\label{d-gue251130yyd}
Let $E: C^\infty(X)\To C^\infty(X)$ be a continuous operator. We say that $E$ is a classical pseudodifferential operator on $X$ of order $m$ if $E$ is smoothing away the diagonal and for any orbifold chart $(U_i, V_i, \Gamma_i, \Psi_i)$, there is a pseudodifferential operator $\tilde E\in L^m_{{\rm cl\,}}(V_i)$ such that for all $f\in C^\infty_0(U_i)$, we have 
\[E(f)(x)=\frac{1}{|\Gamma_i|}\sum_{g\in\Gamma_i}\tilde E(\tilde f)(g\cdot\tilde x)\ \ \mbox{on $U_i$},\]
where $\pi(\tilde x)=x$, $\pi: V_i\rightarrow U_i$ is the natural projection and 
$\tilde f\in C^\infty_0(V_i)$ is the lifting of $f$. 

As in the smooth case, we let $L^m_{{\rm cl\,}}(X)$ denote the space of classical pseudodifferential operators on $X$ of order $m$. 
\end{definition}

Now, we assume that 
\[\square_b: {\rm Dom\,}\square_b\subset L^2(X)\rightarrow L^2(X)\]
has closed range. Let 
\[N: L^2(X)\rightarrow {\rm Dom\,}\square_b\]
be the partial inverse of $\square_b$, that is, 
\begin{equation}\label{e-gue251206yyda}
\begin{split}
&\square_bN+\Pi=I\ \ \mbox{on $L^2(X)$},\\
&N\square_b+\Pi=I\ \ \mbox{on ${\rm Dom\,}\square_b$}.
\end{split}
\end{equation}

\begin{lemma}\label{l-gue251206yyd}
Suppose  that \[\square_b: {\rm Dom\,}\square_b\subset L^2(X)\rightarrow L^2(X)\]
has closed range and the partial inverse $N$ is continuous:
\[N: C^\infty(X)\rightarrow C^\infty(X).\]
Let $(U, V, \Gamma, \Psi)$
be an orbifold chart of $X$. Assume that $X$ is strictly pseudoconvex on $U$. Then, for any $\chi\in C^\infty_0(U)$, $\tau\in C^\infty(X)$ with 
${\rm supp\,}\chi\cap{\rm supp\,}\tau=\emptyset$, we have 
\begin{equation}\label{e-gue251209yyds}\begin{split}
    &\chi\Pi\tau\equiv0\ \ \mbox{on $X$},\\
    &\tau\Pi\chi\equiv0\ \ \mbox{on $X$}.
\end{split}
\end{equation}
\end{lemma} 

\begin{proof}
Let $\hat A$ and $\hat S$ be properly supported continuous operators:
    \[\hat A, \hat S: C^\infty(U)\To C^\infty(U)\]
    given by \eqref{e-gue251206yyd}. We have
    \begin{equation}\label{e-gue251209yyd}
			\begin{split}
	&\mbox{$\square_b\hat{A}+\hat{S}=I$ on $U$},\\
    &\mbox{$\hat{A}^*\square_b+\hat{S}^*=I$ on $U$},\\
				&\square_b\hat{S}\equiv0\ \ \mbox{on $U$},\\
				&\hat{S}\equiv \hat{S}^*\equiv \hat{S}^2\ \ \mbox{on $U$},
			\end{split}
		\end{equation}
		where $\hat{S}^*$ and $\hat{A}^*$ are the formal adjoints of $\hat{S}$ and $\hat{A}$
		with respect to $(\,\cdot\,|\,\cdot\,)$ respectively.  
        Let $\chi\in C^\infty_0(U)$. 
From the second equation of \eqref{e-gue251209yyd} and notice that $\hat A$ and $\hat S$ are properly supported, we have 
\begin{equation}\label{e-gue251209yydI} \mbox{$\chi\hat{S}^*\Pi=\chi\Pi$ on $X$}
\end{equation}
and hence
\begin{equation}\label{e-gue251209yydII}
\mbox{$\Pi\hat{S}\chi=\Pi\chi$ on $X$}.\end{equation}

On the other hand, we have 
\begin{equation}\label{e-gue251209yydIII}\hat S\chi=N\square_b\hat S\chi+\Pi\hat S\chi.\end{equation}
Since $N$ maps smooth functions to smooth functions,
\[N\square_b\hat S\chi: \mcd'(X)\To C^\infty(X)\]
is continuous. Hence 
\[N\square_b\hat S\chi\equiv0.\]
From this observation and \eqref{e-gue251209yydIII}, we get 
\begin{equation}\label{e-gue251209yyda}\hat S\chi=\Pi\hat S\chi+F\ \ \mbox{on $X$},\end{equation}
where $F\equiv0$ on $X$. 

From \eqref{e-gue251209yydII} and \eqref{e-gue251209yyda}, we get 
\[\hat S\chi=\Pi\chi+F\ \ \mbox{on $X$}\]
and hence 
\begin{equation}\label{e-gue251209yydb}
\tau\hat S\chi=\tau\Pi\chi+\tau F,
\end{equation}
where $\tau\in C^\infty(X)$ with 
${\rm supp\,}\chi\cap{\rm supp\,}\tau=\emptyset$. From \eqref{e-gue251209yydb} and notice that $\hat S$ is smoothing away the diagonal, the lemma follows. 
\end{proof}

\begin{remark}\label{r-gue251209yyd}
If Kohn's condition $Y(1)$ holds, then \[\square_b: {\rm Dom\,}\square_b\subset L^2(X)\rightarrow L^2(X)\]
has closed range and the partial inverse $N$ is continuous:
\[N: C^\infty(X)\rightarrow C^\infty(X).\]
\end{remark}

Let $E\in L^m_{{\rm cl\,}}(X)$, $m\leq0$. Consider the Toeplitz operator
\[T_E:=\Pi\circ E\circ\Pi: L^2(X)\To L^2(X).\]
Let $T_E\in\mcd'(X\times X)$ be the distribution kernel of $T_E$. 
The main result of this section is the following asymptotic result for $T_E$ on the strictly pseudoconvex part of $X$.

\begin{theorem}\label{t-gue251209yydf}
Suppose  that \[\square_b: {\rm Dom\,}\square_b\subset L^2(X)\rightarrow L^2(X)\]
has closed range and the partial inverse $N$ is continuous:
\[N: C^\infty(X)\rightarrow C^\infty(X).\]
Let $E\in L^m_{{\rm cl\,}}(X)$, $m\leq0$.
Let $(U, V, \Gamma, \Psi)$
be an orbifold chart of $X$. Assume that $X$ is strictly pseudoconvex on $U$. Then, on $U$, 
		\begin{equation}\label{e-gue251130yydaz}
        T_E(x,y)=\frac{1}{|\Gamma|}\sum_{g,h\in\Gamma}\Bigl(\int _0^{+\infty}e^{it\, \phi(h\cdot\tilde{x},g\cdot \tilde{y})}a_E(h\cdot\tilde{x},g\cdot \tilde{y},t)dt+F(h\cdot\tilde x,g\cdot\tilde y)\Bigr)\ \ \mbox{on $U\times U$} \,,\end{equation}
		where $\pi(\widetilde x)=x$, $\pi(\widetilde y)=y$, $\pi: V\To U$ is the natural projection, $F$ is a smoothing operator on $V$ and 
		\[a_E( \widetilde{x}, \widetilde{y},t)\sim\sum^\infty_{j=0}a_{j,E}( \widetilde{x},  \widetilde{y})\,t^{n+m-j}\]
		in $S^{n+m}_{1, 0}(V\times V\times\mathbb{R}_+)$,
        $a_{j,E}( \widetilde{x},  \widetilde{y})\in C^\infty(V\times V)$. $j=0,1,\ldots$, 
        \[a_{0,E}(\widetilde{x},  \widetilde{x})=\frac{1}{2}\pi^{-n-1}|\det\mathcal{L}_{\widetilde{x}}|e_0(\widetilde x,-\omega_0(\widetilde x)),\ \ \widetilde{x}\in V,\]
        where $\det\mathcal{L}_{\widetilde{x}}=\mu_1(\widetilde x)\cdots\mu_n(\widetilde x)$, $\mu_j(\widetilde x)$, $j=1,\ldots,n$, are the eigenvalues of 
        $\mathcal{L}_{\widetilde{x}}$ with respect to $\langle\,\cdot\,|\,\cdot\,\rangle$, $e_0$ is the principal symbol of $E$, 
		 the phase function $\phi$ is the same as in Theorem~\ref{bdms_thm}. 
         
         Moreover, assume that on $U$, $\langle\,\cdot\,|\,\cdot\,\rangle$ is the Levi metric and the phase $\phi$ satisfies $\mct_{\widetilde y}^2\phi(\widetilde x,\widetilde x) = 0$ for all $\widetilde x \in V$, and satisfies (\ref{phi_props}) (this is always possible). Then we can take $a_{0,E}$ and $a_{1,E}$ so that 
         \begin{equation}\label{e-gue251210ycdt}
         \begin{split}
         &a_{0,E}(\widetilde{x},  \widetilde{x})=\frac{1}{2}\pi^{-n-1}e_0(\widetilde x,-\omega_0(\widetilde x)),\ \ \widetilde{x}\in V,\\
         &\mct_{\widetilde y}a_{0,E}(\widetilde x,\widetilde x)=0,\ \ \widetilde{x}\in V,
         \end{split}
         \end{equation}
         and $a_{1,E}(\widetilde x,\widetilde x)$ satisfies \eqref{main_formula}, for all $\widetilde{x}\in V$. Moreover, if $a_{0,E}$ satisfies \eqref{e-gue251210ycdt}, then $a_{1,E}(\widetilde x,\widetilde x)$ is uniquely determined.
\end{theorem} 

\begin{proof}
Let $\chi, \chi_1, \chi_2\in C^\infty_0(U)$, $\chi_1\equiv1$ on ${\rm supp\,}\chi$, 
$\chi_2\equiv1$ on ${\rm supp\,}\chi_1$.
From Lemma~\ref{l-gue251206yyd}, we
have 
\begin{equation}\label{e-gue251210yyd}
\begin{split}
&\chi_1T_E\chi\\
&=\chi_1\Pi\circ E\circ\Pi\chi\\
&\equiv \chi_1\Pi\chi_2\circ E\circ\chi_1\Pi\chi.
\end{split}
\end{equation}
From the proof of Theorem~\ref{t-gue251130yyd} and \eqref{e-gue251210yyd}, we have 
\begin{equation}\label{e-gue251210yydI}
\begin{split}
&\chi_1T_E\chi\\
&\equiv \chi_1\hat S\chi_2\circ E\circ\chi_1\hat S\chi,
\end{split}
\end{equation}
where $\hat S$ is as in \eqref{e-gue251206yyd}. From \eqref{e-gue251210yydI} and by using complex stationary phase formula of Melin-Sj\"ostrand, we get \eqref{e-gue251130yydaz}.

Since $(\hat S)^2\equiv\hat S$, we can repeat the proof of Theorem~\ref{main_thm} and get the final assertion of the theorem.
\end{proof}

\subsubsection*{Acknowledgements}
We would like to thank the anonymous reviewer for important comments and suggestions.

O. Shabtai would like to express his gratitude to Academia Sinica and National Taiwan University for their hospitality and financial support during visits to Taipei in which parts of this work were conducted.\\
\subsubsection*{Disclosure statement}
The authors report there are no competing interests to declare.
\end{document}